\documentclass[letter, journal,final]{IEEEtran}  

\IEEEoverridecommandlockouts                              

\usepackage{ifpdf}
\usepackage{cite}
\usepackage{amsmath}
\usepackage{amssymb}
\usepackage{amsfonts}
\usepackage{amsthm}
\usepackage{url}
\usepackage{subfig}
\usepackage{fixltx2e}
\usepackage{graphics} 
\usepackage{graphicx}
\usepackage{epsfig} 
\usepackage{stfloats}
\usepackage{caption}
\usepackage{epstopdf}
\usepackage{algorithm,algorithmic}
\usepackage{hyperref}
\hypersetup{
	colorlinks=true,
	linkcolor=blue,
	filecolor=blue,
	urlcolor=blue,
	citecolor = blue,
}
\allowdisplaybreaks

\newtheorem{thm}{Theorem}
\newtheorem{lem}{Lemma}
\newtheorem{rem}{Remark}
\newtheorem*{pf}{Proof}
\newtheorem{prop}{Proposition}
\newtheorem{assum}{Assumption}
\newtheorem{defn}{Definition}

\newcommand{\tr}{\mathrm{tr}}
\newcommand{\T}{^{\top}}
\newcommand{\Ad}[1]{\mathrm{Ad}_{#1}}

\newcommand{\dom}{\text{dom }}
\newcommand{\ie}{\textit{i.e.}}
\newcommand{\g}{\mathsf{g}}

\title{\LARGE \bf
	Hybrid Nonlinear Observers for Inertial Navigation Using Landmark Measurements
}

\author{Miaomiao Wang and Abdelhamid Tayebi 
	\thanks{This work was supported by the National Sciences and Engineering Research Council of Canada (NSERC). A preliminary and partial version of this work was presented in \cite{wang2018navigation}.}
	\thanks{M. Wang is with the Department of Electrical and Computer Engineering, Western University, London, ON N6A 3K7, Canada (e-mail: {\tt\small mwang448@uwo.ca}).}%
	\thanks{A. Tayebi is with the Department of Electrical and Computer Engineering, Western University, London, ON N6A 3K7, Canada, and also with the Department of Electrical Engineering, Lakehead University, Thunder Bay, ON P7B 5E1, Canada
		(e-mail: {\tt\small atayebi@lakeheadu.ca}).}%
}

\begin{document}

\maketitle
\thispagestyle{empty}
\pagestyle{empty}

\begin{abstract}
	This paper considers the problem of attitude, position and linear velocity estimation for rigid body systems relying on inertial measurement unit  and landmark measurements. We propose two hybrid nonlinear observers on the matrix Lie group  $SE_2(3)$, leading to global exponential stability. The first observer relies on fixed gains, while the second one uses variable gains depending on the solution of a continuous Riccati equation. These observers are then extended to handle biased angular velocity and linear acceleration measurements. Both simulation and experimental results are presented to illustrate the performance of the proposed observers.

\end{abstract}

\section{INTRODUCTION}
The development of reliable attitude, position and linear velocity estimation algorithms for inertial navigation systems is instrumental in many applications, such as autonomous underwater and ground vehicles, and unmanned aerial vehicles. It is well known that the attitude of a rigid body can be estimated using body-frame observations of some known inertial vectors obtained, for instance, from an inertial measurement unit (IMU) equipped with a gyroscope, an accelerometer and a magnetometer \cite{bonnabel2008symmetry,mahony2008nonlinear,hua2014implementation}. The position and linear velocity can be obtained, for instance, from a Global Positioning System (GPS)\cite{barczyk2013invariant,barrau2014invariant,grip2013nonlinear}. As a matter of fact, IMU-based nonlinear attitude observers rely on the fact that the accelerometer provides a measurement of the gravity vector in the body-fixed frame, which is true only in the case of negligible linear accelerations.  In applications involving accelerated rigid body systems, a typical solution consists in using linear velocity measurements together with IMU measurements with the so-called velocity-aided attitude observers \cite{bonnabel2009invariant,roberts2011,hua2016stability,hua2017riccati,berkane2017attitude}. However, these observers are not easy to implement in GPS-denied environments (\textit{e.g.,} indoor applications), where it is challenging to obtain the linear velocity.

On the other hand, the pose (position and orientation) can be obtained using estimators relying on inertial-vision systems consisting of an IMU and a stereo vision system attached to the rigid body \cite{rehbinder2003pose}. Most of the existing pose estimators are typically filters of the Kalman-type. Recently, nonlinear geometric observers on Lie groups have made their appearance in the literature \cite{bonnabel2009non,lageman2010gradient}. Nonlinear pose observers designed on $SE(3)$, using landmark and group velocity measurements, have been considered in \cite{vasconcelos2010nonlinear,hua2011observer,hua2015gradient,khosravian2015observers}. However, these observers are shown to guarantee almost global asymptotic stability (AGAS), \ie, the pose error converges to zero from almost all initial conditions except from a set of Lebesgue measure zero. Nonlinear hybrid observers evolving on $SE(3)$  with global asymptotic and exponential stability guarantees have been proposed in \cite{wang2017globally} and \cite{wang2019hybrid}, respectively.

In practice, it is difficult to obtain the linear velocity, using low-cost sensors, in GPS-denied environments. Therefore, developing estimation algorithms that provide the attitude, position and linear velocity, with strong stability guarantees, is of great importance (from theoretical and practical point of views) for inertial navigation systems. It is important to point out that the dynamics of the attitude, position and linear velocity  are not (right or left) invariant, and hence, the extension of the existing invariant observers designed on $SE(3)$ to the estimation problem considered in this work is not trivial. Most of the existing results in the literature, for the state estimation problem for inertial navigation systems, are filters of the Kalman-type, see for instance \cite{mourikis2009vision,panahandeh2014vision}. Recently, an invariant extended Kalman filter (IEKF), using a geometric error on matrix Lie groups, has been proposed in \cite{barrau2017invariant}, and  a Riccati-based geometric pose, linear velocity and gravity direction observer has been proposed in \cite{hua2018riccati}. However, both results are only shown to be locally stable relying on linearizations.

In the present work, we formulate the estimation problem for inertial navigation systems using the matrix Lie group $SE_2(3)$ introduced in \cite{barrau2017invariant}. Then, we propose nonlinear geometric hybrid observers, relying on IMU and landmark measurements, for the estimation of the pose, linear velocity, gyro-bias and accelerometer-bias. The main contributions of this paper can be summarized as follows:
1)  The observers proposed in this paper do not rely on group velocity (angular and linear velocities) measurements as in \cite{vasconcelos2010nonlinear,hua2011observer,hua2015gradient,khosravian2015observers,wang2017globally,wang2019hybrid}.
2) The proposed observers are endowed with global exponential stability guarantees. To the best of our knowledge, there are no results in the literature achieving such strong stability properties for the estimation problem at hand.
3) Contrary to the dynamics on the Lie groups $SO(3), SE(3)$,  the dynamics of the attitude, position and linear velocity on $SE_2(3)$ are not invariant. As a consequence, the application of the hybrid observers proposed in \cite{berkane2017hybrid,wang2019hybrid} to the problem considered in the present paper is not trivial.
4) In our preliminary work \cite{wang2018navigation}, we proposed a fixed-gain hybrid observer for inertial navigation systems in the bias-free case. In this paper, we provide a comprehensive hybrid observer design methodology for inertial navigation systems, with fixed and variable gains, using biased gyro and accelerometer measurements. Moreover, experimental results using IMU and stereo camera measurements are presented to illustrate the performance of the proposed observers.

The rest of this paper is organized as follows: Section \ref{sec:preliminary} introduces some preliminary notions that will be used throughout this paper. Section \ref{sec:observer} is devoted to the design of the hybrid nonlinear observers in the bias-free case. These results are extended to address the problem of biased angular velocity in Section \ref{sec:observerbias}, and biased IMU measurements in Section \ref{sec:biased_IMU}. Simulation and experimental results are presented in Section \ref{sec:simulation} and Section \ref{sec:experimental}, respectively.

\section{Preliminary Material}\label{sec:preliminary}
\subsection{Notations}
The sets of real, non-negative real and natural numbers are denoted as $\mathbb{R}$, $\mathbb{R}^+$ and $\mathbb{N}$, respectively. We denote by $\mathbb{R}^n$ the $n$-dimensional Euclidean space, and denote by $\mathbb{S}^n$ the set of $(n+1)$-dimensional unit vectors. Given two matrices, $A,B\in \mathbb{R}^{m\times n}$, their Euclidean inner product is defined as $\langle\langle A,B\rangle\rangle = \tr(A\T B)$. The Euclidean norm of a vector $x\in \mathbb{R}^n$ is defined as $\|x\| = \sqrt{x\T x}$, and the Frobenius norm of a matrix $X\in \mathbb{R}^{n\times m}$ is given by $\|X\|_F = \sqrt{\langle \langle X, X\rangle\rangle}$. The $n$-by-$n$ identity matrix is denoted by $I_n$. For each $A\in \mathbb{R}^{n\times n}$, we define $\mathcal{E}(A)$ as the set of all eigenvectors of $A$, and $\mathbb{E}(A)\subseteq \mathcal{E}(A)$ as the eigenbasis set of $A$. Let $\lambda^A_i$ be the $i$-th eigenvalue of $A$, and $\lambda^A_{m}$ and $\lambda^A_{M}$ be the minimum and maximum eigenvalue of $A$, respectively. By $\text{blkdiag}(\cdot)$, we denote the block diagonal matrix.

Let $\{\mathcal{I}\}$ be an inertial frame and $\{\mathcal{B}\}$ be a frame attached to a rigid body moving in 3-dimensional space. Define $R\in SO(3)$ as the rotation of frame $\{\mathcal{B}\}$ with respect to frame $\{\mathcal{I}\}$, where $SO(3):=\{R\in \mathbb{R}^{3\times 3}| RR\T =R\T R= I_3, \det(R)=1\}$. Let $p\in \mathbb{R}^3$ and $v\in \mathbb{R}^3$ be the position and linear velocity of the rigid-body expressed in the inertial frame $\{\mathcal{I}\}$, respectively. We consider the following extended Special Euclidean group of order 3 proposed in \cite{barrau2017invariant}:   $SE_{2}(3) := SO(3) \times  \mathbb{R}^3 \times \mathbb{R}^3  \subset \mathbb{R}^{5\times 5}$, which is defined as
\begin{equation}
SE_{2}(3) = \{X= \mathcal{T}(R,v,p) |  R\in SO(3), p, v\in \mathbb{R}^3 \}
\end{equation}
where the map $\mathcal{T}: SO(3)\times \mathbb{R}^3 \times \mathbb{R}^3 \to SE_2(3)$ is defined by  \cite{barrau2017invariant}
\begin{align*}
\mathcal{T}(R,v,p) &=\begin{bmatrix}
R & v & p\\
0_{1\times 3} & 1 & 0 \\
0_{1\times 3} & 0 & 1
\end{bmatrix}.
\end{align*}
For every $X=\mathcal{T}(R,v,p)$, one has $X^{-1} = \mathcal{T}(R\T,-R\T v, -R\T p)$. Denote $T_X SE_2(3)\in \mathbb{R}^{5\times 5}$ as the \textit{tangent space} of $SE_2(3)$ at point $X$.
The \textit{Lie algebra} of  $SE_2(3)$, denoted by $\mathfrak{se}_2(3)$, is given by
\begin{align*}
\mathfrak{se}_2(3) :=& \left\{ \left. U= \begin{bmatrix}
\Omega  & \alpha & v\\
0_{2\times 3}  & 0_{2\times 3} & 0_{2\times 3}
\end{bmatrix}\in \mathbb{R}^{5\times 5} \right|  \right.\\
&\qquad \qquad \qquad \qquad \qquad \left.     \Omega \in \mathfrak{so}(3), v,\alpha\in \mathbb{R}^3\right\}
\end{align*}
where $\mathfrak{so}(3) = \{\Omega \in \mathbb{R}^{3\times 3}| \Omega = -\Omega\T \}$ denoting the Lie algebra of $SO(3)$.
Let $\times$ be the vector cross-product on $\mathbb{R}^3$ and define the map $(\cdot)^\times: \mathbb{R}^3 \to \mathfrak{so}(3)$ such that $x\times y = x^\times y, \forall x,y\in \mathbb{R}^3 $.
Define the inverse isomorphism of $(\cdot)^\times$ as $\text{vec}: \mathfrak{so}(3)\to \mathbb{R}^3$ such that $\text{vec}(\omega^\times)=\omega, (\text{vec}(\Omega))^\times = \Omega , \forall \omega\in \mathbb{R}^3, \Omega \in \mathfrak{so}(3)$.  For any $R\in SO(3)$, we define $|R|_I\in [0,1]$ as the normalized Euclidean distance on $SO(3)$ with respect to the identity $I_3$, such that $|R|_I^2 = \frac{1}{8}\|I_3-R\|_F^2=\frac{1}{4}\tr(I_3-R)$. Let the map $\mathcal{R}_a: \mathbb{R}\times \mathbb{S}^2$ represent the well-known angle-axis parametrization of the attitude, which is given by $\mathcal{R}_a(\theta,u):=I_3 + \sin\theta u^\times + (1-\cos\theta) (u^\times)^2$ with $\theta$ being the rotation angle and $u$ the rotation axis.
For any matrix $A_1\in \mathbb{R}^{3\times 3}$, define $\mathbb{P}_a(A_1)$ as the anti-symmetric projection of $A_1$, such that $\mathbb{P}_a(A_1) = (A_1-A_1\T)/2$. For a matrix $A_1 \in \mathbb{R}^{3\times 3}$, we define $\psi(A_1) = \text{vec}(\mathbb{P}_a(A_1))$. Then, one has the identity $  \langle\langle A_1, u^\times\rangle\rangle = 2 u\T \psi(A_1)$. Let  $\mathbb{P}:\mathbb{R}^{5\times 5} \rightarrow \mathfrak{se}_2(3)$ denote the projection of $A$ on the Lie algebra  $\mathfrak{se}_2(3)$, such that, for all $U\in \mathfrak{se}_2(3)$, $ A \in \mathbb{R}^{5\times 5}$ one has
$ \langle\langle A, U\rangle\rangle=\langle\langle U,\mathbb{P}( A)\rangle\rangle = \langle\langle \mathbb{P}(A),U\rangle\rangle
$.
For all $A_1\in \mathbb{R}^{3\times 3}, a_2, \cdots, a_5 \in
\mathbb{R}^3$ and $a_6,\cdots,a_9\in \mathbb{R}$, one has
\begin{equation}
\mathbb{P}\left(\begin{bmatrix}
A_1 & a_2 & a_3 \\
a_4\T & a_6 & a_7 \\
a_5\T & a_8 & a_9
\end{bmatrix}\right) = \begin{bmatrix}
\mathbb{P}_a(A_1) & a_2 & a_3\\
0_{1\times 3}  & 0 & 0 \\
0_{1\times 3}  & 0 & 0
\end{bmatrix}. \label{eqn:project_P}
\end{equation}
Given a rigid body with configuration $X\in SE_2(3)$, for all $X\in SE_2(3), U\in \mathfrak{se}_2(3)$, the \textit{adjoint map} $\Ad{}: SE_2(3) \times \mathfrak{se}_2(3) \to \mathfrak{se}_2(3)$ is given by $ \Ad{X}U   := X U X^{-1}$. For all $X_1,X_2\in SE_2(3), U\in \mathfrak{se}_2(3)$, one can verify that $\Ad{X_1}\Ad{X_2} U= \Ad{X_1X_2} U$.

\subsection{Hybrid Systems Framework}
Define the \textit{hybrid time domain} as the subset $E \subset \mathbb{R}^+ \times \mathbb{N}$ in the form of
$
E = \bigcup_{j=0}^{J-1} ([t_j,t_{j+1}] \times \{j\}),
$
for some finite time sequence $0=t_0 \leq t_1 \leq \cdots \leq t_J$, with the last interval possibly in the form $([t_{J},T)\times \{J\})$ with $T$ finite or $T=\infty$. On each hybrid time domain there is a natural ordering of points : $(t,j)\preceq (t',j')$ if $t\leq t'$ and $j\leq j'$.
Given a manifold $\mathcal{M}$ embedded in $\mathbb{R}^n$ and its tangent space $T\mathcal{M}$, we consider the following hybrid system \cite{goebel2009hybrid}:
\begin{equation}\mathcal{H}:  ~~~
\begin{cases}
\dot{x} ~~= F(x),& \quad x \in \mathcal{F}   \\
x^{+} = G(x),& \quad x \in \mathcal{J}
\end{cases} \label{eqn:hybrid_system}
\end{equation}
where the \textit{flow map} $F: \mathcal{M} \to T\mathcal{M}$ describes the continuous flow of $x$ on the \textit{flow set} $\mathcal{F} \subset \mathcal{M}$; the \textit{jump map} $G: \mathcal{M} \rightrightarrows \mathcal{M}$ describes the discrete jump of $x$ on the \textit{jump set} $\mathcal{J} \subset \mathcal{M}$. A hybrid arc is a function $x: \dom{x} \to \mathcal{M}$, where $\dom{x}$ is a hybrid time domain and, for each fixed $j$, $t \mapsto x(t,j)$ is a locally absolutely continuous function on the interval $I_j = \{t:(t,j) \in \dom{x}\}$.
Note that $x^+$  denotes the value $x$ after a jump, namely, $x^+=x(t,j+1)$ with $x(t,j)$ denoting the value of $x$ before the jump.
For more details on dynamic hybrid systems, we refer the reader to \cite{goebel2009hybrid,goebel2012hybrid} and references therein. For the hybrid system $\mathcal{H}$, a closed set $\mathcal{A}\subset \mathcal{M}$ is said to be (locally) exponentially stable if there exist strictly positive scalars $\kappa, \lambda$ and $\mu$ such that, for any initial condition $|x(0,0)|_{\mathcal{A}} < \mu$, each solution $x$ to $\mathcal{H}$ satisfies $|x(t,j)|_{\mathcal{A}} \leq \kappa \exp(-\lambda(t+j))|x(0,0)|_{\mathcal{A}}$ for all $(t,j)\in \dom x$, and it is said to be globally exponentially stable if one allows $\mu \to + \infty$  \cite{teel2013lyapunov}.

\subsection{Continuous Riccati Equation}
Consider  $A(t)\in  \mathbb{R}^{n\times n} $ and $ C(t)\in \mathbb{R}^{m\times n}$ as matrix-valued functions of time $t$, and  suppose that $A(t)$ and $C(t)$   are continuous and bounded for all $t\geq 0$. Consider the following continuous Riccati equation (CRE):
\begin{equation}
\dot{P} = A(t)P + PA(t)\T - PC(t)\T Q(t) C(t)P + V(t) \label{eqn:CRE}
\end{equation}
where $P(0) \in \mathbb{R}^{n\times n}$ is a symmetric positive definite matrix and $V(t)\in \mathbb{R}^{n\times n}$ and $Q(t)\in \mathbb{R}^{m\times m}$ are  uniformly positive definite matrices.   The following definition formulates the well-known uniform observability condition in terms of the Gramian matrix:
\begin{defn}\label{def:UCO}
	The pair $(A(t),C(t))$  is uniformly   observable  if there exist constants $\delta, \mu >0$ such that $\forall t\geq 0 $
	\begin{equation}
	W(t,t+\delta)  :=\frac{1}{\delta} \int_{t}^{t+\delta} \Phi(\tau,t)\T C(\tau)\T   C(\tau) \Phi(\tau,t) d\tau \geq \mu I_n    \label{eqn:gramiancondition}
	\end{equation}
	with $\Phi(t,\tau)$ being  the \textit{state transition matrix} associated with $A(t)$
\end{defn}
 Note that     $W(t,t+\delta)$ is naturally upper bounded by some constant since that matrices $A(t)$ and $C(t)$ are bounded for all $t\geq 0$ by assumption. To establish global existence, uniqueness and boundedness of the solution  of the CRE (\ref{eqn:CRE}), sufficient conditions are presented in the following lemma:
\begin{lem}[\cite{bucy1967global,bucy1972riccati}] \label{lemma:gramian}
	If there exist constants $\delta,\mu_q, \mu_v>0$ such that $\forall t\geq 0$
	\begin{align}
	&\frac{1}{\delta} \int_{t}^{t+\delta} \Phi(t+\delta,\tau) V(\tau) \Phi(t+\delta,\tau)\T d\tau      \geq \mu_v I_n   \label{eqn:gramianV} \\
	&\frac{1}{\delta} \int_{t}^{t+\delta} \Phi(\tau,t)\T C(\tau)\T Q(\tau) C(\tau) \Phi(\tau,t) d\tau  \geq \mu_q I_n     \label{eqn:gramianQ}
	\end{align}
	then the solution $P(t)$ of the CRE  \eqref{eqn:CRE} is well-defined on $\mathbb{R}^+$, and there exist   constants $0< p_m\leq p_M < \infty$ such that  $p_m I_n \leq P(t) \leq p_M I_n $.
\end{lem}
Note that conditions (\ref{eqn:gramianV}) and (\ref{eqn:gramianQ}) are satisfied if   $V(t)\in \mathbb{R}^{n\times n}$ and $Q(t)\in \mathbb{R}^{m\times m}$ are  uniformly positive definite and the pair $(A(t),C(t))$ is uniformly   observable.
In the traditional Kalman filter, $V(t)$ and $Q^{-1}(t)$  are interpreted as   covariance matrices of additive noise on the system state and output.

\subsection{Kinematics and Measurements}
Consider the following kinematics of a rigid body navigating in three-dimensional space:
\begin{align}
\dot{R} & = R \omega^\times  \label{eqn:R}\\
\dot{p} & = v  \label{eqn:p}\\
\dot{v} & = \mathsf{g} + Ra  \label{eqn:v}
\end{align}
where $\mathsf{g}\in \mathbb{R}^3$ denotes the gravity vector,  $\omega\in \mathbb{R}^3$ denotes the angular velocity expressed in the body-fixed frame, and $a\in \mathbb{R}^3$ is the body-frame ``apparent acceleration" capturing all non-gravitational forces applied to the rigid body expressed in the body-fixed frame. We assume that $\omega$ and $a$ are continuous and available for measurement. In this paper, we consider the configuration of the rigid body represented by an element of the matrix Lie group $X=\mathcal{T}(R,v,p) \in SE_2(3)$. Let us introduce the nonlinear map $f: SE_2(3) \times \mathbb{R}^3 \times \mathbb{R}^3 \to T_XSE_2(3)$, such that the kinematics (\ref{eqn:R})-(\ref{eqn:v}) can be rewritten in the following compact form:
\begin{equation}
\dot{X} = f(X,\omega,a):=\begin{bmatrix}
R\omega^\times& \mathsf{g}+ Ra   & v \\
0_{1\times 3} & 0 & 0\\
0_{1\times 3} & 0 & 0
\end{bmatrix}. \label{eqn:dX}
\end{equation}

Consider a family of $n$ landmarks available for measurement, and let $p_i \in \mathbb{R}^3$ be the position of the $i$-th landmark expressed in frame $\{\mathcal{I}\}$. The landmark measurements expressed in the   frame $\{\mathcal{B}\}$ are denoted as
\begin{equation}
y_i := R\T(p_i-p), \quad i=1,2,\cdots,n. \label{eqn:output_y}
\end{equation}
The three-dimensional landmark position measurements can be obtained, for instance, using bearing measurements generated from a stereo vision system \cite[Eq. (26)]{hamel2018riccati}.
Let $r_i := [p_i\T~0~1]\T \in \mathbb{R}^{5}, i=1,\cdots,n$ be the new inertial reference vectors expressed in the   frame $\{\mathcal{I}\}$, and $b_i := [y_i\T ~0~1]\T\in \mathbb{R}^{5}, i=1,\cdots,n$ be their measurements expressed in the  frame $\{\mathcal{B}\}$. From (\ref{eqn:output_y}), one has
\begin{equation}
b_i  = h(X,r_i) :=X^{-1} r_i, \quad i=1,2,\cdots,n. \label{eqn:output_X}
\end{equation}
Note that, the Lie group action $h:SE_2(3)\times \mathbb{R}^{5}\to \mathbb{R}^{5}$ is a \textit{right group action} in the sense that for all $X_1,X_2\in SE_2(3)$ and $r\in \mathbb{R}^5$, one has
$
h(X_2,h(X_1,r) ) = h(X_1X_2,r) .
$
For later use, we define $ r:= [r_1 ~r_2~\cdots ~r_n]\in \mathbb{R}^{5\times n}$ and $b :=  [b_1 ~b_2~\cdots ~b_n]\in \mathbb{R}^{5\times n}$.

\begin{assum}\label{assum:1}
	Assume that there exist at least three non-collinear landmarks among the $n\geq 3$ measurable landmarks.
\end{assum}
Assumption \ref{assum:1} is common in pose estimation on $SE(3)$ using landmark measurements  \cite{hua2011observer,vasconcelos2010nonlinear,hua2015gradient,khosravian2015observers,wang2017globally,wang2019hybrid}.
Define the matrix $M:=\textstyle \sum_{i=1}^n  k_i (p_i - p_c)(p_i - p_c)\T$ with $k_i>0, \forall i = 1,2,\cdots,n$, $k_c := \sum_{i=1}^n k_i$ and  $p_c :=   \frac{1}{k_c}\sum_{i=1}^{n} k_i p_i$. The matrix $M$ can be rewritten as $M =A_1 - k_cp_cp_c\T$ with $A_1 := \sum_{i=1}^n k_i p_ip_i\T$.
Given three non-collinear landmarks, it is always possible to guarantee that the matrix $M$ is positive semi-definite with no more than one zero eigenvalue through an appropriate choice of the gains $k_i, i=1,2\cdots, n$.

\section{Hybrid Observers Design Using Bias-free Angular Velocity}\label{sec:observer}
\subsection{Continuous Observer and Undesired Equilibria}
Let $\hat{X}:=\mathcal{T}(\hat{R},\hat{v},\hat{p})\in SE_2(3)$ be the estimate of the state $X$, where $\hat{R}$ denotes the estimate of the attitude $R$, $\hat{v}$ denotes the estimate of the linear velocity $v$ and $\hat{p}$ denotes the estimate of the position $p$. Define the right-invariant estimation error as $\tilde{X} := X\hat{X}^{-1} = \mathcal{T}(\tilde{R},\tilde{v},\tilde{p})$ with $\tilde{R}:=R\hat{R}\T, \tilde{v}: = v - \tilde{R}\hat{v}$ and $\tilde{p}: = p- \tilde{R}\hat{p}$.
Consider the following time-invariant continuous observer:
\begin{align} 	
&\dot{\hat{X}} = f(\hat{X},\omega,a) - \Delta \hat{X}  \label{eqn:observer_smooth}\\
&\Delta :=
-\Ad{X_c} \left(\mathbb{P}(X_c^{-1}(r-\hat{X}b) K_n r\T   X_c^{-\top} K)  \right)
\label{eqn:innovation_term}	
\end{align}
where $\hat{X}(0)\in SE_2(3)$ and $X_c : = \mathcal{T}(I_3, 0_{3\times 1} ,p_c) \in SE_2(3)$ with $ p_c$ given before. The gain parameters are given by
\begin{equation}\small 
K_n = \text{diag}(k_1,\cdots,k_n),~~   K = \begin{bmatrix}
k_R I_3 & 0_{3\times 1} & 0_{3\times 1}\\
0_{1\times 3} & 0 & 0\\
0_{1\times 3} & k_v & k_p
\end{bmatrix}    \label{eqn:def_Kn_K}
\end{equation}
with constant scalars $k_R,k_p,k_v, k_i >0, i=1,2,\cdots,n$.

\begin{rem}
	Note that the proposed continuous observer is designed on the matrix Lie group $SE_2(3)$ directly, which is different from most of the existing Kalman-type filters. The observer has two parts: the term $f(\hat{X},\omega,a)$ relying on the measurements of $\omega$ and $a$, and an innovation term $\Delta$ designed in terms of the estimated state $\hat{X}$ and landmarks measurements.
\end{rem}	
\begin{rem}
A homogeneous transformation matrix $X_c \in SE_2(3)$ is introduced in the innovation term $\Delta$, which intends to transform the inertial reference vectors to a specific frame. Considering the transformation $\bar{r}_i= X_c^{-1} r_i,  i = 1,2,\cdots,n$, the innovation term $\Delta$ defined in (\ref{eqn:innovation_term}) can be simplified as
$\Delta=-\Ad{X_c}  (\mathbb{P}((\bar{r}-X_c^{-1}\hat{X}b) K_n \bar{r}\T K)  )$ with $\bar{r}=[\bar{r}_1,\cdots,\bar{r}_n]$. Choosing $X_c = \mathcal{T}(I_3, 0_{3\times 1} ,p_c)$,   leads to a nice decoupling property in the closed loop dynamics, which will be discussed later. Similar techniques can be found in \cite{wang2017globally,wang2019hybrid}.	
\end{rem}	

Let $\tilde{y}_i: = p_i - \hat{p} - \hat{R}y_i = (I_3-\tilde{R}\T)p_i + \tilde{R}\T \tilde{p}  $ for all $i=1,2,\cdots,n$. From the definitions of $r, b$ and $K_n$, one obtains
\begin{align*}
& X_c^{-1}(r-\hat{X}b)K_n r\T  X_c^{-\top}
 \\
&\qquad \qquad  = \begin{bmatrix}
\sum_{i=1}^n k_i \tilde{y}_i (p_i-p_c)\T  & 0_{3\times 1} & \sum_{i=1}^n k_i \tilde{y}_i\\
0_{1\times 3} & 0 & 0\\
0_{1\times 3} & 0 & 0
\end{bmatrix}
\end{align*}
where we made use of the fact $(r-\hat{X}b)K_n r\T  = \sum_{i=1}^n k_i(r_i-\hat{X}b_i)r_i\T$.
Then, from the definitions of $K$, the Adjoint map $\Ad{}$ and the projection map $\mathbb{P}$, the term $\Delta$ defined in (\ref{eqn:innovation_term}) becomes
\begin{align}
\Delta
&= -\begin{bmatrix}
k_R \mathbb{P}_a( \Delta_R ) & k_v \Delta_p & k_p \Delta_p -k_R \mathbb{P}_a( \Delta_R)p_c\\
0_{1\times 3} & 0 & 0\\
0_{1\times 3} & 0 & 0
\end{bmatrix}   \label{eqn:Delta}
\end{align}
with
\begin{align}
\Delta_R & := \sum_{i=1}^n k_i \tilde{y}_i (p_i - p_c)\T
= (I_3 -\tilde{R})\T M \label{eqn:Delta_R}\\
\Delta_p  &  :=  \sum_{i=1}^n k_i\tilde{y}_i
 = k_c\tilde{R}\T (\tilde{p} - (I_3-\tilde{R})p_c)  \label{eqn:Delta_p}
\end{align}
where we made use of the facts  $\tilde{y}_i     = (I_3-\tilde{R}\T)p_i + \tilde{R}\T \tilde{p}$, $ k_c = \sum_{i=1}^n k_i, \sum_{i=1}^n k_i  p_i = k_c p_c$ and $M = \sum_{i=1}^n k_i  p_ip_i\T - k_c p_c p_c\T$. From (\ref{eqn:Delta_R}) one verifies that $\Delta_{R}$  is strictly bounded. It is important to mention that $\Delta_R   =  (I_3-\tilde{R}\T)\sum_{i=1}^n k_i  p_i p_i\T  + k_c \tilde{R}\T \tilde{p} p_c\T $ when choosing $X_c = I_5$. This implies that the position estimation error $\tilde{p}$ will affect the attitude estimation error $\tilde{R}$ when the weighted center of landmarks  is not located at the origin (\ie, $p_c\neq 0_{3\times 1}$). By choosing $X_c = \mathcal{T}(I_3, 0_{3\times 1} ,p_c)$, one can show that $\Delta_R$ in (\ref{eqn:Delta_R}) is independent from the position estimation error $\tilde{p}$.

For the sake of simplicity, let us define the new position estimation error $\tilde{p}_e: = \tilde{p} - (I-\tilde{R})p_c$. In view of (\ref{eqn:dX}), (\ref{eqn:observer_smooth}) and (\ref{eqn:Delta}), one has the following closed-loop system:
\begin{equation}
\begin{cases}
\dot{\tilde{R}}  ~= \tilde{R} (-k_R\mathbb{P}_a(M\tilde{R}))  \\
\dot{\tilde{p}}_e =    -k_pk_c\tilde{p}_e + \tilde{v} \\
\dot{\tilde{v}}~=  -  k_v k_c  \tilde{p}_e + (I_3-\tilde{R})\mathsf{g}
\end{cases} \label{eqn:closed-loop_F}
\end{equation}
where we made use of the facts   $\mathbb{P}_a( \Delta_R ) = \mathbb{P}_a( (I_3-\tilde{R}\T)M) = \mathbb{P}_a( M\tilde{R} ) $ and $\dot{\tilde{p}}_e = \dot{\tilde{p}} - k_R\tilde{R}\mathbb{P}_a(\Delta_R)p_c$. It is clear that  $\mathcal{T}(\tilde{R},\tilde{v},\tilde{p}) = I_5$ if and only if $\mathcal{T}(\tilde{R},\tilde{v},\tilde{p}_e) = I_5$. Note that the   errors $\tilde{v}$ and $\tilde{p}_e$ considered in this paper are different from the linear errors (\ie, $v-\hat{v}$ and $p-\hat{p}$) considered in the classical EKF-based navigation filters. The modified geometric errors lead to an interesting decoupling property for the closed-loop system, where the dynamics of $\tilde{R}$ are not dependent on $\tilde{p}_e$ and $\tilde{v}$ as shown in the first equation of (\ref{eqn:closed-loop_F}).

\begin{prop} \label{pro:pro_1}
	Consider the closed-loop dynamics (\ref{eqn:closed-loop_F}). Let $\Psi$ be the set of undesired equilibrium points (\ie, all the equilibrium points except $I_5$) of the closed-loop dynamics, which is given by
	\begin{align}
	\Psi :=& \left\{   \mathcal{T}(\tilde{R},\tilde{v},\tilde{p}_e) \in SE_2(3)| \tilde{R}=\mathcal{R}_a(\pi,u), u\in \mathcal{E}(M),\right.    \nonumber \\
	& \qquad \quad  \left.      \tilde{p}_e =   {(k_c k_v)^{-1}} (I-\tilde{R}) \g , \tilde{v}= k_c k_p \tilde{p}_e \right\} . \label{eqn:Psi}
	\end{align}
\end{prop}
\begin{pf}
	The proof of   Proposition \ref{pro:pro_1} is straightforward.  From the first equation of (\ref{eqn:closed-loop_F}), identity $\dot{\tilde{R}}=0_{3\times 3}$ implies that $\mathbb{P}_a(M\tilde{R})=0$. Applying \cite[Lemma 2]{mayhew2011synergistic}, one can further show that $\tilde{R}\in  \{I_3\} \cup \Psi_M$, where the set $\Psi_M:= \{\tilde{R}\in SO(3): \tilde{R} = \mathcal{R}_a(\pi,u), u \in \mathcal{E}(M)\}$ denotes the set of undesired equilibrium points of the rotational error dynamics. Substituting $\tilde{R}\in \Psi_M$ into the identities $\dot{\tilde{v}}=0_{3\times 1}$ and $\dot{\tilde{p}}_e=0_{3\times 1}$,  one can easily verify (\ref{eqn:Psi}) from (\ref{eqn:closed-loop_F}).
\end{pf}
\begin{rem}
	From the dynamics of $\tilde{R}$ in (\ref{eqn:closed-loop_F}), it is easy to verify that the equilibrium point $(\tilde{R}=I_3)$ is almost globally asymptotically stable \cite{mahony2008nonlinear}.
	It is important to mention that, due to the topology of the Lie group $SO(3)$ as pointed out in \cite{koditschek1989application}, it is impossible to achieve robust and global stability results with smooth (or even discontinuous) state observers \cite{mayhew2011synergistic,mayhew2013synergistic}. Hence, the best stability result one can achieve with the continuous observer (\ref{eqn:observer_smooth})-(\ref{eqn:innovation_term}) is AGAS. This motivates the design of hybrid observers leading to robust and global stability results as shown in the next section.
\end{rem}

\subsection{Fixed-Gain Hybrid Observer Design}
Define the following real-valued cost function $\Upsilon:SE_2(3) \times \mathbb{R}^{5\times n} \times \mathbb{R}^{5\times n} \to \mathbb{R}^+$
\begin{equation}
\Upsilon(\hat{X},r,b)  :=  \frac{1}{2}\sum_{i=1}^n k_i \|  (r_i-r_c) - \hat{X} (b_i-b_c)\|^2
\label{eqn:Phi_R}
\end{equation}
where $r_c: =  \sum_{i=1}^n \frac{k_i}{k_c} r_i = [p_c\T~ 0~1]\T$ and $b_c: =   \sum_{i=1}^n \frac{k_i}{k_c} b_i = [ y_c\T  ~ 0~1]\T$ with  $ y_c:=\sum_{i=1}^n \frac{k_i}{k_c}y_i = R\T(p_c - p) $. From the definitions of $p_c, k_c$ and $M$, one can rewrite $\Upsilon(\hat{X},r,b)$ as
$
\Upsilon(\hat{X},r,b)   =  \frac{1}{2}\sum_{i=1}^n k_i \| (p_i-p_c) - \hat{R}(y_i-y_c)  \|^2
 = \tr((I_3 - \tilde{R})M).
$
Given a non-empty and finite \textit{transformation set} $\mathbb{Q}\subset SE_2(3)$, let us define the real-valued function $\mu_\mathbb{Q}: SE_2(3) \times \mathbb{R}^{5\times n} \times \mathbb{R}^{5\times n} \to \mathbb{R}$ as
\begin{equation}
\mu_\mathbb{Q}(\hat{X},r,b): = \Upsilon(\hat{X},r,b) - \min_{X_q\in \mathbb{Q}} \Upsilon(X_q^{-1}\hat{X},r,b). \label{eqn: definition_mu_R}
\end{equation}
The flow set $\mathcal{F}_o$ and jump set $\mathcal{J}_o$ are defined as follows:
\begin{align}
\mathcal{F}_o&: = \{\hat{X} \in SE_2(3) ~  |~ \mu_\mathbb{Q}(\hat{X},r,b)  \leq \delta  \}  \label{eqn:F_map}\\
\mathcal{J}_o&: = \{\hat{X} \in SE_2(3)  ~ |~ \mu_\mathbb{Q}(\hat{X},r,b)  \geq \delta  \}  \label{eqn:J_map}
\end{align}
with a constant scalar $\delta>0$ and a non-empty finite set $\mathbb{Q}$ to be designed later. One can show that the sets $\mathcal{F}_o$ and $\mathcal{J}_o$ are closed, and $\mathcal{F}_o\cup \mathcal{J}_o= SE_2(3)$.  We propose the following hybrid observer:
\begin{align}
&\mathcal{H}^{o}_1: \begin{cases}
\dot{\hat{X}} = f(\hat{X},\omega,a) - \Delta \hat{X}, & \hat{X}\in \mathcal{F}_o   \\
\hat{X}^+  = X_q^{-1} \hat{X},~~  X_q \in \gamma(\hat{X}),  & \hat{X}\in \mathcal{J}_o
\end{cases}
\label{eqn:observer_X} \\
&\Delta :=
-\Ad{X_c}  (\mathbb{P}(X_c^{-1}(r-\hat{X}b) K_n r\T   X_c^{-\top} K)  )
\label{eqn:innovation_term2}
\end{align}
where $\hat{X}(0)\in SE_2(3)$ and the set-valued map $\gamma: SE_2(3) \rightrightarrows SE_2(3)$ is defined by
\begin{equation}
\gamma(\hat{X}) := \left\{X_q\in \mathbb{Q}| X_q = \arg \min_{X_q \in \mathbb{Q}} \Upsilon (X_q^{-1} \hat{X},r,b) \right\}   \label{eqn:gamma}.
\end{equation}

We define the extended space and state as $\mathcal{S}^c_1: = SE_2(3) \times SO(3) \times \mathbb{R}^3 \times \mathbb{R}^3 \times \mathbb{R}^+$ and $x^c_1:=(\hat{X},\tilde{R},\tilde{p}_e,\tilde{v},t)$, respectively.
In view of  (\ref{eqn:closed-loop_F}) and (\ref{eqn:set Q})-(\ref{eqn:gamma}), one obtains the following hybrid closed-loop system:
\begin{equation}
\mathcal{H}^c_1:  \begin{cases}
\dot{x}^c_1 ~~= F_1(x^c_1), &x^c_1\in \mathcal{F}^c_1  \\
{x^c_1}^+ = G_1(x^c_1),  & x^c_1 \in \mathcal{J}^c_1
\end{cases} \label{eqn:closed-loop}
\end{equation}
with $\mathcal{F}^c_1:=\{x^c_1=(\hat{X},\tilde{R},\tilde{p}_e,\tilde{v},t)\in \mathcal{S}^c_1: \hat{X}\in \mathcal{F}_o\}$, $\mathcal{J}^c_1:=\{x^c_1=(\hat{X},\tilde{R},\tilde{p}_e,\tilde{v},t)\in \mathcal{S}^c_1: \hat{X}\in \mathcal{J}_o\}$, and
\begin{equation*}
F_1(x^c_1)= \begin{pmatrix}
f(\hat{X},\omega,a) - \Delta\hat{X}\\
\tilde{R} (-k_R\mathbb{P}_a(M\tilde{R}))\\
-k_c k_p\tilde{p}_e + \tilde{v}   \\
-k_c k_v \tilde{p}_e  + (I-\tilde{R})\mathsf{g}\\
1
\end{pmatrix}, G_1(x^c_1) = \begin{pmatrix}
X_q^{-1} \hat{X} \\
\tilde{R}R_q \\
\tilde{p}_e \\
\tilde{v} \\
t
\end{pmatrix}
\end{equation*}
where $X_q = \mathcal{T}(R_q,v_q,p_q)\in \gamma(\hat{X})$, and we made use of the facts: $\tilde{R}^+ =\tilde{R}R_q$, $\tilde{p}^+ = \tilde{p}+\tilde{R}(I_3-R_q)p_c, \tilde{p}_e^+ =\tilde{p}_e$ and $\tilde{v}^+ = \tilde{v}$.
Note that the sets $\mathcal{F}^c_1$ and $\mathcal{J}^c_1$ are closed, and $\mathcal{F}^c_1 \cup \mathcal{J}^c_1 = \mathcal{S}^c_1$. Note also that the hybrid system $\mathcal{H}^c_1$   satisfies the hybrid basic conditions of \cite{goebel2009hybrid} and is autonomous by taking $\omega$ and $a$ as functions of time $t$.

The main idea behind our hybrid observer is the introduction of a resetting mechanism  to avoid all the undesired equilibrium points of the closed-loop system (\ref{eqn:closed-loop}) in the flow set $\mathcal{F}^c_1$, \textit{i.e.}, all the undesired equilibrium points of the closed-loop system lie in the jump set $\mathcal{J}^c_1$. The innovation term $\Delta$ and the transformation set $\mathbb{Q}$ are designed to guarantee a decrease of a Lyapunov function in both flow set $\mathcal{F}^c_1$ and jump set $\mathcal{J}^c_1$. The transformation set $\mathbb{Q}$ is given by
\begin{align}
\mathbb{Q} := & \{X = \mathcal{T}(R,v,p) \in SE_2(3) | R=\mathcal{R}_a(\theta,u)  , u\in \mathbb{U} , \nonumber \\
& \qquad \qquad \qquad \qquad  ~~~ p = (I_3-R)p_c, v=0_{3\times 1}\}   \label{eqn:set Q}
\end{align}	
with a constant $\theta \in (0,\pi]$ and  a finite set $\mathbb{U} \subset \mathbb{S}^2$, which can be chosen as per one of the following two approaches:
\begin{itemize}
	\item [1)] A superset of  the eigenbasis set of $M$, \textit{i.e.}, $\mathbb{U}\supseteq \mathbb{E}(M)$, if $\lambda_1^M \geq  \lambda_2^M \geq \lambda_3^M >0$ or $\lambda_1^M >  \lambda_2^M > \lambda_3^M = 0$.
	\item [2)] A set that contains any three orthogonal unit vectors in $\mathbb{R}^3$, if   $\tr(M) -2\lambda_{\max}^M> 0$.
\end{itemize}
Note that   approach 1) considers the case where the matrix $M$ is positive definite    or positive semi-definite with distinct eigenvalues, however it requires the information about the eigenvectors of $M$. Note also that approach 2) does not need any information about the eigenvectors of $M$, but it requires a strong condition $\tr(M) -2\lambda_{\max}^M> 0$.

\begin{lem}[\cite{wang2019hybrid}]\label{lemma:Delta_M}
	Let $M=M\T$ be a positive semi-definite matrix under Assumption \ref{assum:1} and $\mathbb{U}$ be a nonempty finite set of unit vectors. Consider the map $\Delta_M: \mathbb{S}^2 \times \mathbb{U} \to \mathbb{R}$ defined as
	\begin{equation}
	\Delta_M(u,v) :=   u\T (\tr(M_v)I_3- M_v) u  \label{eqn:definition_Delta_uv}
	\end{equation}
	where  $M_v := M(I_3-2vv\T)$ and $v\in \mathcal{E}(M)$.  Define the constant scalar
	\begin{equation}
	\Delta_M^*:= \min_{v\in \mathcal{E}(M)} \max_{u\in \mathbb{U}} \Delta_M(u,v) . \label{eqn:Delta_M_star}
	\end{equation}
	Then, the following results hold:	
	\begin{itemize}
		\item [1)] Let $\mathbb{U}$ be a superset of $\mathbb{E}(M)$ (\textit{i.e.}, $\mathbb{U}\supseteq \mathbb{E}(M)$), then the following inequality holds:
		\begin{align*}
		\Delta_M^*  \geq \begin{cases}
		\frac{2}{3}  \lambda_1^M & \textit{if } \lambda_1^M = \lambda_2^M = \lambda_3^M >0\\
		\min \{2\lambda^M_{1} , \lambda_3^M  \} &  \textit{if } \lambda_1^M = \lambda_2^M \neq \lambda_3^M >0 \\
		\tr(M) - \lambda_{\max}^M  &  \textit{if } \lambda_1^M > \lambda_2^M > \lambda_3^M \geq 0
		\end{cases}. \label{eqb:solution_delta_U}
		\end{align*}
		\item [2)] Let $M$ be a matrix such that $\tr(M) -2\lambda_{\max}^M> 0$, and let $\mathbb{U}$ be a set that contains any three orthogonal unit vectors in $\mathbb{R}^3$, then the following inequality holds:
		\begin{equation*}
		\Delta_M^*  \geq  {  \textstyle  \frac{2}{3}} (\tr(M)  - 2\lambda_{\max}^M  ) .\label{eqn:Delta_M2}
		\end{equation*}	
	\end{itemize}	
\end{lem}

\begin{prop}\label{pro:undesired_eq}
	Consider the hybrid system $\mathcal{H}^c_1$ defined in (\ref{eqn:closed-loop}) with the set $\mathbb{Q}$ defined in (\ref{eqn:set Q}). Suppose that Assumption \ref{assum:1} holds.
	Then, for all $\delta<(1-\cos\theta) \Delta^*_M $ with $\theta$ given in (\ref{eqn:set Q}) and $\Delta_M^*$ given as per Lemma \ref{lemma:Delta_M}, one has $SE_2(3) \times\Psi \times \mathbb{R}^+ \subseteq \mathcal{J}^c_1$.
\end{prop}

See Appendix \ref{sec:undesired_eq} for the proof. Proposition \ref{pro:undesired_eq} provides a choice for the gap $\delta$, ensuring that the set of undesired equilibrium points of the flow dynamics of (\ref{eqn:closed-loop}) is a subset of the jump set $\mathcal{J}^c_1$, \ie, all the undesired equilibrium points of the flows of $\mathcal{H}^c_1$ lie in the jump set $\mathcal{J}^c_1$.

Let us define the closed set $\mathcal{A}_1:=\{ (\hat{X},\tilde{R},\tilde{p}_e, \tilde{v},t)\in \mathcal{S}^c_1:  \tilde{R}=I_3,\tilde{p}_e=\tilde{v}=0_{3\times 1}\}$.
Now, one can state one of our main results.
\begin{thm}\label{theo:theo_1}
	Consider the hybrid system $\mathcal{H}^c_1$ defined in (\ref{eqn:closed-loop}). Suppose that Assumption \ref{assum:1} holds. For the sets $\mathcal{F}_o$ and $\mathcal{J}_o$, choose the set $\mathbb{Q}$ as in (\ref{eqn:set Q}) and $\delta<(1-\cos(\theta))\Delta_M^*$ with $\theta$ given in (\ref{eqn:set Q}) and $\Delta_M^*$ given as per Lemma \ref{lemma:Delta_M}. Then, the number of discrete jumps is finite and the set $\mathcal{A}_1$ is uniformly globally exponentially stable.		
\end{thm}
\begin{pf}
	See Appendix \ref{sec:theo_1}
\end{pf}


\begin{rem}
	In view of (\ref{eqn:Delta}), (\ref{eqn:observer_X}), (\ref{eqn:gamma}) and (\ref{eqn:set Q}), the proposed hybrid observer can be explicitly expressed, in terms of the available measurements, as follows:
	\begin{align*}
	&\left.\begin{array}{ll}
	\dot{\hat{R}} &= \hat{R}\omega^\times + k_R \mathbb{P}_a( \Delta_R) \hat{R} \\
	\dot{\hat{p}} &= \hat{v} +  k_R \mathbb{P}_a( \Delta_R) (\hat{p} -p_c) +  k_p \Delta_p     \\
	\dot{\hat{v}} &= \mathsf{g}+ \hat{R}a +  k_R \mathbb{P}_a( \Delta_R) \hat{v} + k_v \Delta_p
	\end{array}\right\}  \hat{X}\in \mathcal{F}_o   \\
	&\left.\begin{array}{l}
	\hat{R}^+  = R_q\T \hat{R}    \\
	\hat{p}^+  ~= R_q\T (\hat{p} - (I_3-R_q)p_c)\\
	\hat{v}^+  ~= R_q\T \hat{v}
	\end{array}     \right\}   \hat{X}\in \mathcal{J}_o
	\label{eqn:observer_Rpv}
	\end{align*}
	where $R_q= \min_{R_q\in \mathcal{R}_a(\theta,\mathbb{U}) } \sum_{i=1}^n k_i \| (p_i-p_c) - R_q\T \hat{R}(y_i -y_c)  \|^2$ with  $\mathcal{R}_a(\theta,\mathbb{U})=\{R\in SO(3): R = \mathcal{R}_a(\theta,u), u\in \mathbb{U} \}$,  $\Delta_R$ and $\Delta_p$  are defined in (\ref{eqn:Delta_R}) and (\ref{eqn:Delta_p}). From (\ref{eqn:set Q}), it is always possible to  guarantee the existence of $R_q$, which belongs to $\mathcal{R}_a(\theta,\mathbb{U})$ with $\theta$ and $\mathbb{U}$ given in (\ref{eqn:set Q}).
\end{rem}

\subsection{Variable-Gain Hybrid Observer Design}\label{sec:V_ghod}

In this subsection we provide a different version of the hybrid observer $\mathcal{H}^o_1$ using  variable gains relying on the solution of a CRE.
Let us define the following gain map $\mathbb{P}_{\mathcal{K}}: \mathbb{R}^{5\times 5} \to \mathfrak{se}_2(3)$ inspired by \cite{khosravian2015observers}, such that for all $A_1\in \mathbb{R}^{3\times 3}, a_2, \cdots, a_5 \in
\mathbb{R}^3$ and $a_6,\cdots,a_9\in \mathbb{R}$, one has
\begin{equation}\small
\mathbb{P}_{\mathcal{K}}\left(\begin{bmatrix}
A_1 & a_2 & a_3 \\
a_4\T & a_6 & a_7 \\
a_5\T & a_8 & a_9
\end{bmatrix}\right) = \begin{bmatrix}
k_R\mathbb{P}_a(A_1) & K_v a_2 & K_p a_3\\
0_{1\times 3}  & 0 & 0 \\
0_{1\times 3}  & 0 & 0
\end{bmatrix}  \label{eqn:project_PK}
\end{equation}
where $\mathcal{K} := (k_R, K_p, K_v)$ with $k_R>0$ and $K_p, K_v \in \mathbb{R}^{3\times 3}$ to be designed. Then, we propose the following hybrid observer
\begin{align}
&\mathcal{H}^{o}_2: \begin{cases}
 \begin{array}{l}
\dot{\hat{X}} = f(\hat{X},\omega ,a) - \Delta \hat{X}
\end{array} & \hat{X}\in \mathcal{F}_o   \\
 \begin{array}{l}
\hat{X}^+  = X_q^{-1} \hat{X}, \quad X_q \in \gamma(\hat{X})
\end{array}   & \hat{X}\in \mathcal{J}_o
\end{cases}
\label{eqn:observer_X2} \\
&~ \Delta :=
-\Ad{X_c}  (\mathbb{P}_{\mathcal{K}}(X_c^{-1}(r-\hat{X}b) K_n r\T   X_c^{-\top})  )
\label{eqn:innovation_term4}	
\end{align}
where $\hat{X}(0)\in SE_2(3)$. The map $\gamma $ is defined in (\ref{eqn:gamma}), and the flow and jump sets $\mathcal{F}_o$ and $\mathcal{J}_o$ are defined in (\ref{eqn:F_map}) and (\ref{eqn:J_map}), respectively. Note that the main difference between the hybrid observers $\mathcal{H}^o_1$ and $\mathcal{H}^o_2$ is the innovation term $\Delta$. Instead of using constant scalar gains $k_v, k_p$ as in the observer $\mathcal{H}^o_1$, the new observer $\mathcal{H}^o_2$ uses variable matrix gains $K_v, K_p$ to be designed later via a CRE.

In view of (\ref{eqn:dX}), and (\ref{eqn:project_PK})-(\ref{eqn:innovation_term4}), one has the following closed-loop system in the flows:
\begin{equation}
\begin{cases}
\begin{array}{ll}
\dot{\tilde{R}} &= \tilde{R} (-k_R\mathbb{P}_a(M\tilde{R}))  \\
\dot{\tilde{p}}_e &=    -k_c\tilde{R}K_p \tilde{R}\T \tilde{p}_e + \tilde{v} \\
\dot{\tilde{v}}&=  -  k_c \tilde{R} K_v  \tilde{R}\T \tilde{p}_e + (I_3-\tilde{R})\mathsf{g}
\end{array}
\end{cases} \label{eqn:closed-loop_F2}.
\end{equation}
Define the new variable $\mathsf{x}:= [(R\T\tilde{p}_e)\T, (R\T\tilde{v})\T ]\in \mathbb{R}^6$. It is easy to verify that $\|\tilde{p}_e\|=\|\tilde{v}\|=0$ if $\|\mathsf{x}\|=0$. Let $L := k_c[\hat{R}\T K_p\T \hat{R}, \hat{R}\T K_v\T \hat{R}]\T \in \mathbb{R}^{6\times 3}$ and $\nu = [0_{1\times 3} ~	\mathsf{g}\T(R -\hat{R})]\T \in \mathbb{R}^{6}$. From (\ref{eqn:closed-loop_F2}) one obtains the dynamics of $\mathsf{x}$ as
\begin{equation}
\dot{\mathsf{x}} = A(t) \mathsf{x} - L C \mathsf{x} + \nu \label{eqn:dx}
\end{equation}
with
\begin{equation}
A(t) := \begin{bmatrix}
-\omega(t)^\times & I_3 \\
0_{3\times 3} & -\omega(t)^\times
\end{bmatrix},\quad  C  := \begin{bmatrix}
I_3 &
 0_{3\times 3}
\end{bmatrix} . \label{eqn:pair_A_C}
\end{equation}
Note that the dynamics of $\mathsf{x}$ are linear time-varying and $\nu$ can be viewed as a disturbance term. The variable-gain matrix $L$ can be updated as $L = PC\T Q(t)$, with $P$ being the solution of the  CRE (\ref{eqn:CRE}). 	Let $L_1,L_2\in \mathbb{R}^{3\times 3}$ such that $ [L_1\T, L_2\T]\T=L $. Then, the gain matrices $K_p$ and $K_v$ can be computed as
	\begin{equation}
	  K_p =\frac{1}{k_c}  \hat{R} L_1 \hat{R}\T,\quad  K_v = \frac{1}{k_c}  \hat{R} L_2 \hat{R}\T \label{eqn:Kpv2}.
	\end{equation}
	For the sake of simplicity, one can always choose the weights $k_i>0, i=1,\cdots,n$ such that $k_c = \sum_{i=1}^n k_i = 1$.

\begin{lem}\label{lem:existenceP}
	The pair $(A(t),C)$ defined in (\ref{eqn:pair_A_C}) is uniformly   observable.
\end{lem}
See Appendix \ref{sec:existenceP} for the proof.
%

Define the extended space $\mathcal{S}^c_2: = SE_2(3)  \times SO(3) \times \mathbb{R}^6  \times \mathbb{R}^+$ and the extended state $x^c_2:=(\hat{X},\tilde{R},\mathsf{x},t)$. Let us define the set $\mathcal{A}_2:=\{ (\hat{X},\tilde{R},\mathsf{x},t)\in \mathcal{S}^c_2:    \tilde{R}=I_3, \mathsf{x}=0_{6\times 1} \}$.  Now, one can state the following result:
\begin{thm}\label{theo:theo_2}
	Consider the inertial navigation system (\ref{eqn:dX})-(\ref{eqn:output_y}) with the hybrid observer (\ref{eqn:observer_X2})-(\ref{eqn:innovation_term4}). Suppose that Assumption \ref{assum:1} holds.  For the sets $\mathcal{F}_o$ and $\mathcal{J}_o$, choose the set $\mathbb{Q}$ as in (\ref{eqn:set Q}) and $\delta<(1-\cos(\theta))\Delta_M^*$ with $\theta$ given in (\ref{eqn:set Q}) and $\Delta_M^*$ given as per Lemma \ref{lemma:Delta_M}. Let $k_R>0$, $Q(t)$ and $V(t)$ be  strictly positive definite. Then, the number of discrete jumps is finite and the set $\mathcal{A}_2$ is uniformly globally exponentially stable.		
\end{thm}
\begin{pf}
	See Appendix \ref{sec:theo_2}.
\end{pf}

\section{Hybrid Observers Design Using Biased Angular Velocity}\label{sec:observerbias}
\subsection{Fixed-Gain Hybrid Observer Design}
In the previous section, nonlinear hybrid observers have been designed using non-biased angular velocity measurements.
In this section, we consider the case where the angular velocity measurements contain an unknown constant or slowly varying bias. Let $b_\omega$ be the constant unknown angular velocity bias, such that $\omega_y = \omega + b_\omega$. Define $\hat{b}_\omega$ as the estimate of $b_\omega$ and $\tilde{b}_\omega:=\hat{b}_\omega-b_\omega$ as the estimation error.

We propose the following hybrid nonlinear observer for inertial navigation with biased angular velocity:
\begin{align}
& \mathcal{H}^{o}_3:
\begin{cases}\left.\begin{array}{l}
\dot{\hat{X}} = f(\hat{X},\omega_y-\hat{b}_\omega ,a) - \Delta \hat{X} \\
\dot{\hat{b}}_\omega  = -k_\omega \hat{R}\T \psi(\Delta_R)
\end{array}\right\}   (\hat{X},\hat{b}_\omega)\in \mathcal{F}_o \times \mathbb{R}^3
\\
\left.\begin{array}{l}
\hat{X}^+  = X_q^{-1} \hat{X},\quad X_q \in \gamma(\hat{X}) \\
\hat{b}_\omega^+~ = \hat{b}_\omega
\end{array}\right\}  ~ (\hat{X},\hat{b}_\omega)\in \mathcal{J}_o \times \mathbb{R}^3
\end{cases}
\label{eqn:observer_X3} \\
&\Delta :=
-\Ad{X_c}  (\mathbb{P}(X_c^{-1}(r-\hat{X}b) K_n r\T X_c^{-\top}  K )  )
\label{eqn:innovation_term3}
\end{align}
where $\hat{X}(0) \in SE_2(3), \hat{b}_\omega(0) \in \mathbb{R}^3$, $k_\omega>0$, $K_n$ and $K$ are given by (\ref{eqn:def_Kn_K}) and $\Delta_R$ is given in (\ref{eqn:Delta_R}). The map $\gamma $ is defined in (\ref{eqn:gamma}) and the flow and jump sets $\mathcal{F}_o, \mathcal{J}_o$ are defined in (\ref{eqn:F_map}) and (\ref{eqn:J_map}), respectively.

Consider the extended space and state as $\mathcal{S}^c_3: = \mathcal{S}^c_1  \times \mathbb{R}^3   \times \mathbb{R}^3 $ and $x^c_3:=(x^c_1, \hat{b}_\omega ,\tilde{b}_\omega)$. Let us define the   set $\mathcal{A}_3:=\{ (x^c_1,\hat{b}_\omega, \tilde{b}_\omega )\in \mathcal{S}^c_3: x^c_1 \in \mathcal{A}_1,  \tilde{b}_\omega = 0_{3\times 1} \}$. Let $|x^c_3|_{\mathcal{A}_3}\geq 0$ denote the distance to the set $\mathcal{A}_3$ such that
$
|x^c_3|_{\mathcal{A}_3}^2  : = \textstyle \inf_{y=(\bar{X},I_3,0,0,\bar{t},\bar{b}_\omega,0)\in \mathcal{A}_3}  (\|\bar{X}-\hat{X}\|_F^2 + |\tilde{R}|_I^2
+ \|\tilde{p}_e\|^2+ \|\tilde{v}\|^2 +\|\bar{t}-t\|^2+\|\bar{b}_\omega - \hat{b}_\omega\|^2 +\|\tilde{b}_\omega\|^2  )
= \textstyle |\tilde{R}|_I^2   + \|\tilde{p}_e\|^2 + \|\tilde{v}\|^2 + \|\tilde{b}_\omega\|^2 .
$

Before stating our next result, the following assumption is made:
\begin{assum}\label{assum:2}
	The state $X$ and angular velocity $\omega$ are uniformly bounded.
\end{assum}

\begin{thm}\label{theo:theo_3}
	Consider the inertial navigation system (\ref{eqn:dX})-(\ref{eqn:output_y}) with the hybrid observer (\ref{eqn:observer_X3})-(\ref{eqn:innovation_term3}). Suppose that Assumption \ref{assum:1} and Assumption \ref{assum:2} hold.  For the sets $\mathcal{F}_o$ and $\mathcal{J}_o$, choose the set $\mathbb{Q}$ as in (\ref{eqn:set Q}) and $\delta<(1-\cos(\theta))\Delta_M^*$ with $\theta$ given in (\ref{eqn:set Q}) and $\Delta_M^*$ given as per Lemma \ref{lemma:Delta_M}. Let $k_R, k_p, k_v, k_\omega>0$. Then, for any initial condition $x^c_3(0,0) \in \mathcal{S}^c_3$ the number of discrete jumps is finite, and the solution   $x^c_3(t,j)$ is complete and there exist $\kappa,\lambda_F>0$ (depending on the initial conditions) such that
	\begin{equation}
	|x^c_3(t,j)|^2_{\mathcal{A}_3} \leq \kappa\exp\left(- \lambda_F(t+j)\right) |x^c_3(0,0)|^2_{\mathcal{ A}_3}
	\end{equation}
	for all $(t,j)\in \dom x^c_3$.
\end{thm}
\begin{pf}
	See Appendix \ref{sec:theo_3}.
\end{pf}
\begin{rem}
	Note that the parameters $\lambda_F$ and $\kappa$ depend on the initial conditions, which is different from Theorem \ref{theo:theo_1}. This non-uniform type of exponential stability is a consequence of the angular velocity bias (see, for instance, the hybrid observers on $SO(3)$ in \cite{berkane2017hybrid} and the hybrid observers on $SE(3)$ in \cite{wang2019hybrid}).
\end{rem}

\subsection{Variable-Gain Hybrid Observer Design}
We propose the following Riccati-based hybrid nonlinear observer for inertial navigation with biased angular velocity:
\begin{align}
&\mathcal{H}^{o}_4:
\begin{cases}\left.\begin{array}{l}
\dot{\hat{X}} = f(\hat{X},\omega_y-\hat{b}_\omega ,a) - \Delta \hat{X} \\
\dot{\hat{b}}_\omega = -k_\omega\hat{R}\T \psi(\Delta_R)
\end{array}\right\}  (\hat{X},\hat{b}_\omega)\in \mathcal{F}_o \times \mathbb{R}^3
\\
\left.\begin{array}{l}
\hat{X}^+  = X_q^{-1} \hat{X},\quad X_q \in \gamma(\hat{X}) \\
\hat{b}_\omega^+~ = \hat{b}_\omega
\end{array}\right\} ~ (\hat{X},\hat{b}_\omega)\in \mathcal{J}_o \times \mathbb{R}^3
\end{cases}
\label{eqn:observer_X4} \\
&\Delta :=-\Ad{X_c}  (\mathbb{P}_{\mathcal{K}}  (X_c^{-1}(r-\hat{X}b) K_n r\T  X_c^{-\top}))
\label{eqn:innovation_term5}
\end{align}
where $\hat{X}(0) \in SE_2(3), \hat{b}_\omega(0)  \in \mathbb{R}^3$, $k_\omega>0$, $K_n$ is given by (\ref{eqn:def_Kn_K}) and $\Delta_R$ is given in (\ref{eqn:Delta_R}). The gain map $\mathbb{P}_{\mathcal{K}}$ is given by (\ref{eqn:project_PK}). The map $\gamma $ is defined in (\ref{eqn:gamma}) and the flow and jump sets $\mathcal{F}_o$ and $\mathcal{J}_o$ are defined in (\ref{eqn:F_map}) and (\ref{eqn:J_map}), respectively.

In view of (\ref{eqn:dX}), (\ref{eqn:Delta}), and (\ref{eqn:observer_X4})-(\ref{eqn:innovation_term5}), one has the following closed-loop system in the flows:
\begin{equation*}
\begin{cases}
\begin{array}{l}
\dot{\tilde{R}}~=  \tilde{R} ((\hat{R}\tilde{b}_\omega)^\times - k_R\mathbb{P}_a(M\tilde{R}))  \\
\dot{\tilde{b}}_\omega = -k_\omega \hat{R}\T \psi(M\tilde{R}) \\
\dot{\tilde{p}}_e =    -k_c\tilde{R}K_p \tilde{R}\T \tilde{p}_e + \tilde{v}- (R\tilde{b}_\omega)^\times (p-p_c-\tilde{p}_e)  \\
\dot{\tilde{v}}~=  -  k_c \tilde{R} K_v  \tilde{R}\T \tilde{p}_e + (I_3-\tilde{R})\mathsf{g} - (R\tilde{b}_\omega)^\times (v-\tilde{v})
\end{array}
\end{cases}.
\end{equation*}
From the definition of $\mathsf{x}$ and matrix $L$ in Section \ref{sec:V_ghod}, one has the dynamics of $\mathsf{x}$ as
\begin{equation}
\dot{\mathsf{x}} =   A(t) \mathsf{x} -  L C  \mathsf{x} + \nu \label{eqn:dx2}
\end{equation}
with matrix $C$ defined in (\ref{eqn:pair_A_C}), vector $\nu = [
((R\T (p-p_c))^\times  \tilde{b}_\omega)\T,
((R\T v)^\times \tilde{b}_\omega + (I_3-\tilde{R})\mathsf{g})\T
]\T$,  and
\begin{equation}
A(t)  := \begin{bmatrix}
- (\omega_y-\hat{b}_\omega)^\times   & I_3 \\
0_{3\times 3} &  - (\omega_y-\hat{b}_\omega) ^\times
\end{bmatrix} . \label{eqn:pair_Abar_C}
\end{equation}
The gain matrix $L$ can be updated by $L =   PC\T Q(t)$, with $P$ being the solution of the  CRE (\ref{eqn:CRE}). Note that matrices $K_p$ and $K_v$ can be easily obtained from (\ref{eqn:Kpv2}) with $L=PC\T Q(t)$.

\begin{lem}\label{lem:existenceP2}
	The pair $(A(t),C )$ with $A(t)$ defined in (\ref{eqn:pair_Abar_C}) and $C$ defined in (\ref{eqn:pair_A_C}) is uniformly   observable.
\end{lem}
The proof of Lemma \ref{lem:existenceP2} can be conducted  using similar steps as in the proof of Lemma \ref{lem:existenceP}, by introducing  the matrices
$T(t) = \text{blkdiag}([
\bar{R}(t),  \bar{R}(t),\bar{R}(t)
])$, $ S(t) = \text{blkdiag}([
(-\omega_y(t)+\hat{b}_\omega(t))^\times, (-\omega_y(t)+\hat{b}_\omega(t))^\times,  (-\omega_y(t)+\hat{b}_\omega(t))^\times
])$  and constant matrix $\bar{A} = A(t) - S(t)$ with the  rotation matrix $\bar{R}(t)$   generated by  $\dot{\bar{R}}(t) = (-\omega_y(t)+\hat{b}_\omega(t))^\times \bar{R}(t)  $ and  $\bar{R}(0)\in SO(3)$.

Define the extended space and state as $\mathcal{S}^c_4: = \mathcal{S}^c_2  \times \mathbb{R}^3   \times \mathbb{R}^3 $ and $x^c_4:=(x^c_2, \hat{b}_\omega ,\tilde{b}_\omega)$. Let us define the   set $\mathcal{A}_4:=\{ (x^c_2,\hat{b}_\omega, \tilde{b}_\omega )\in \mathcal{S}^c_4: x^c_2 \in \mathcal{A}_2,  \tilde{b}_\omega = 0_{3\times 1} \}$. Let $|x^c_4|_{\mathcal{A}_4}\geq 0$ denote the distance to the set $\mathcal{A}_4$ such that
$
|x^c_4|_{\mathcal{A}_4}^2  : = \textstyle \inf_{y=(\bar{X},I_3,0,\bar{t},\bar{b}_\omega,0)\in \mathcal{A}_4}  (\|\bar{X}-\hat{X}\|_F^2 + |\tilde{R}|_I^2
+ \|\mathsf{x}\|^2 +\|\bar{t}-t\|^2 +\|\bar{b}_\omega-\hat{b}_\omega\|^2+\|\tilde{b}_\omega\|^2)
= \textstyle |\tilde{R}|_I^2 + \|\mathsf{x}\|^2 + \|\tilde{b}_\omega\|^2.
$
Now, one can state the following result:

\begin{thm}\label{theo:theo_4}
	Consider the hybrid observer (\ref{eqn:observer_X4})-(\ref{eqn:innovation_term5}) for the system (\ref{eqn:dX})-(\ref{eqn:output_y}). Suppose that Assumption \ref{assum:1} and Assumption \ref{assum:2} hold.  For the sets $\mathcal{F}_o$ and $\mathcal{J}_o$, choose the set $\mathbb{Q}$ as in (\ref{eqn:set Q}) and $\delta<(1-\cos(\theta))\Delta_M^*$ with $\theta$ given in (\ref{eqn:set Q}) and $\Delta_M^*$ given as per Lemma \ref{lemma:Delta_M}. Let $k_R, k_\omega>0$,  $Q(t)$ and $ V(t)$  be strictly positive definite. Then, for any initial condition $x^c_4(0,0) \in \mathcal{S}^c_4$ the number of discrete jumps is finite, and  the solution   $x^c_4(t,j)$ is complete and there exist $\kappa, \lambda_F>0$ (depending on the initial conditions) such that
	\begin{equation}
	|x^c_4(t,j)|^2_{\mathcal{A}_4} \leq \kappa\exp\left(- \lambda_F(t+j)\right) |x^c_4(0,0)|^2_{\mathcal{ A}_4}
	\end{equation}
	for all $(t,j)\in \dom x^c_4$.	
\end{thm}

\begin{pf}
	See Appendix \ref{sec:theo4}.
\end{pf}

\section{Hybrid Observer Design Using Biased Angular Velocity and Linear Acceleration}\label{sec:biased_IMU}
In this section, we consider the case where the acceleration measurements contain an unknown constant or slowly varying bias. Let $b_a$ be the unknown acceleration bias  such that $a_y = a +b_a$. Define   $\hat{b}_a$ as the estimation of $b_a$ and $\tilde{b}_a := \hat{b}_a - b_a$ as the estimation error.

We propose the following Riccati-based hybrid nonlinear observer for inertial navigation:
\begin{align}
&\mathcal{H}^{o}_5:
\begin{cases}\underbrace{\begin{array}{lr}
	\dot{\hat{X}} = f(\hat{X},\omega_y-\hat{b}_\omega ,a_y-\hat{b}_a) - \Delta \hat{X} \\
	\dot{\hat{b}}_\omega = -k_\omega\hat{R}\T \psi(\Delta_R) \\
	\dot{\hat{b}}_a = - \hat{R}\T K_a  \Delta_p
	\end{array} }_{(\hat{X},\hat{b}_\omega,\hat{b}_a)\in \mathcal{F}_o \times \mathbb{R}^3 \times \mathbb{R}^3}
\\
\underbrace{\begin{array}{l}
	\hat{X}^+  = X_q^{-1} \hat{X},~~ X_q \in \gamma(\hat{X}) \\
	\hat{b}_\omega^+~ = \hat{b}_\omega\\
	\hat{b}_a^+~ =  \hat{b}_a
	\end{array}}_{ (\hat{X},\hat{b}_\omega,\hat{b}_a)\in \mathcal{J}_o \times \mathbb{R}^3 \times \mathbb{R}^3  }
\end{cases}
\label{eqn:observer_X5} \\
&\Delta :=-\Ad{X_c}  (\mathbb{P}_{\mathcal{K}} (X_c^{-1}(r-\hat{X}b) K_n r\T  X_c^{-\top}))
\label{eqn:innovation_term6}
\end{align}
where $\hat{X}(0) \in SE_2(3), \hat{b}_\omega(0), \hat{b}_a(0)  \in \mathbb{R}^3$, $k_\omega>0$, $K_n$ is given by (\ref{eqn:def_Kn_K}) and the innovation terms $\Delta_R, \Delta_p$ are given in (\ref{eqn:Delta_R}) and (\ref{eqn:Delta_p}), respectively. The gain map $\mathbb{P}_{\mathcal{K}}$ is given by (\ref{eqn:project_PK}). The map $\gamma $ is defined in (\ref{eqn:gamma}) and the flow and jump sets $\mathcal{F}_o, \mathcal{J}_o$ are defined in (\ref{eqn:F_map}) and (\ref{eqn:J_map}), respectively.
In view of (\ref{eqn:dX}), (\ref{eqn:Delta}), and (\ref{eqn:observer_X5})-(\ref{eqn:innovation_term6}), one has the following closed-loop system in the flows:
\begin{equation}\small
\left\{
\begin{array}{l}
\dot{\tilde{R}}~=  \tilde{R} ((\hat{R}\tilde{b}_\omega)^\times - k_R\mathbb{P}_a(M\tilde{R}))\\
\dot{\tilde{b}}_\omega = -k_\omega \hat{R}\T \psi(M\tilde{R})\\
\dot{\tilde{p}}_e =    -k_c\tilde{R}K_p \tilde{R}\T \tilde{p}_e + \tilde{v} - (R\tilde{b}_\omega)^\times (p-p_c-\tilde{p}_e)\\
\dot{\tilde{v}}~=  -  k_c \tilde{R} K_v  \tilde{R}\T \tilde{p}_e + R\tilde{b}_a + (I_3-\tilde{R})\mathsf{g} - (R\tilde{b}_\omega)^\times (v-\tilde{v})\\
\dot{\tilde{b}}_a  =      -k_c  \hat{R}\T  K_a  \tilde{R}\T \tilde{p}_e
\end{array}
\right.
. \label{eqn:closed-loop_F5}
\end{equation}
Define the new variable $ {\mathsf{x}}:=[(R\T \tilde{p}_e)\T, (R\T \tilde{v})\T,\tilde{b}_a\T]\T\in \mathbb{R}^9$. Note that   $ \|{\mathsf{x}}\| = 0$ if and only if $\|\tilde{p}_e\|=\|\tilde{v}\|=\|\tilde{b}_a\|=0 $. In view of (\ref{eqn:closed-loop_F5}),  one has the following dynamics of $ {\mathsf{x}}$:
\begin{equation}
\dot{ {\mathsf{x}}} = A(t)  {\mathsf{x}} -  LC   {\mathsf{x}} +  {\nu} \label{eqn:dx3}
\end{equation}
with $L:= k_c[\hat{R}\T K_p\T \hat{R},\hat{R}\T K_v \hat{R},\hat{R}\T K_a \hat{R}]\T\in \mathbb{R}^{9\times 3}$, $ \nu: = [((R\T (p-p_c))^\times  \tilde{b}_\omega)\T,~
((R\T v)^\times \tilde{b}_\omega + (I_3-\tilde{R})\mathsf{g})\T,0_{1\times 3}]\T \in \mathbb{R}^9$, and
\begin{align}
A(t) & := \begin{bmatrix}
- (\omega_y-\hat{b}_\omega)^\times & I_3 & 0_{3\times 3} \\
0_{3\times 3} &  - (\omega_y-\hat{b}_\omega)^\times & I_3 \\
  0_{3\times 3} & 0_{3\times 3} & 0_{3\times 3}
\end{bmatrix} \nonumber \\
C &:=  \begin{bmatrix}
I_3~ 0_{3\times 3} ~ 0_{3\times 3}
\end{bmatrix}. \label{eqn:pair_Abar_C2}
\end{align}
The gain matrix $L$ can be updated by $L =   PC\T Q(t)$, with $P$ being the solution of the  CRE (\ref{eqn:CRE}).
Let $L_1,L_2, L_3\in \mathbb{R}^{3\times 3}$ such that $ [L_1\T, L_2\T, L_3\T]\T =L $. Then, the gain matrices $K_p$, $K_v$ and $K_a$ can be computed as
\begin{equation*}
K_p =\frac{1}{k_c}  \hat{R} L_1 \hat{R}\T,  K_v = \frac{1}{k_c}  \hat{R} L_2 \hat{R}\T,   K_a = \frac{1}{k_c}  \hat{R} L_3 \hat{R}\T.
\end{equation*}

\begin{assum}\label{assum:3}
	The time-derivative of $\omega$ is uniformly bounded.
\end{assum}

\begin{lem}\label{lem:existenceP3}
	The pair $(A(t),C )$ defined in (\ref{eqn:pair_Abar_C2}) is uniformly   observable under Assumption \ref{assum:2} and Assumption \ref{assum:3}.
\end{lem}
See Appendix \ref{sec:existenceP3} for the proof.

Define the extended space and state as $\mathcal{S}^c_5: = SE_2(3)  \times SO(3) \times \mathbb{R}^9  \times \mathbb{R}^+   \times \mathbb{R}^3   \times \mathbb{R}^3  \times \mathbb{R}^3   \times \mathbb{R}^3 $ and $x^c_5:=(\hat{X}, \tilde{R}, {\mathsf{x}},t, \hat{b}_\omega ,\tilde{b}_\omega,\hat{b}_a,\tilde{b}_a)$. Let us define the   set $\mathcal{A}_5:=\{ (\hat{X}, \tilde{R}, {\mathsf{x}},t, \hat{b}_\omega ,\tilde{b}_\omega,\hat{b}_a,\tilde{b}_a)\in \mathcal{S}^c_5: \tilde{R}=I_3,  {\mathsf{x}}=0_{9\times 1},  \tilde{b}_\omega = \tilde{b}_a=0_{3\times 1} \}$. Let $|x^c_5|_{\mathcal{A}_5}\geq 0$ denote the distance to the set $\mathcal{A}_5$ such that
$
|x^c_5|_{\mathcal{A}_5}^2  : = \textstyle \inf_{y=(\bar{X},I_3,0,\bar{t},\bar{b}_\omega,0,\bar{b}_a,0)\in \mathcal{A}_5}  (\|\bar{X}-\hat{X}\|_F^2 + |\tilde{R}|_I^2
+ \| {\mathsf{x}}\|^2 +\|\bar{t}-t\|^2 +\|\bar{b}_\omega-\hat{b}_\omega\|^2+\|\tilde{b}_\omega\|^2+\|\bar{b}_a-\hat{b}_a\|^2+\|\tilde{b}_a\|^2)
= \textstyle |\tilde{R}|_I^2 + \|\mathsf{x}\|^2 + \|\tilde{b}_\omega\|^2 + \|\tilde{b}_a\|^2.
$

Now, one can state the following result:
\begin{thm}\label{theo:theo_5}
	Consider the hybrid observer (\ref{eqn:observer_X5})-(\ref{eqn:innovation_term6}) for the system (\ref{eqn:dX})-(\ref{eqn:output_y}). Suppose that Assumption \ref{assum:1} - Assumption \ref{assum:3} hold.  For the sets $\mathcal{F}_o$ and $\mathcal{J}_o$, choose the set $\mathbb{Q}$ as in (\ref{eqn:set Q}) and $\delta<(1-\cos(\theta))\Delta_M^*$ with $\theta$ given in (\ref{eqn:set Q}) and $\Delta_M^*$ given as per Lemma \ref{lemma:Delta_M}. Let $k_R, k_\omega>0$ and $Q(t)$, $V(t)$ be strictly positive definite. Then, for any initial condition $x^c_5(0,0) \in \mathcal{S}^c_5$ the number of discrete jumps is finite, and  the solution   $x^c_5(t,j)$ is complete and there exist $\kappa, \lambda_F>0$ (depending on the initial conditions) such that
	\begin{equation}
	|x^c_5(t,j)|^2_{\mathcal{A}_5} \leq \kappa\exp\left(- \lambda_F(t+j)\right) |x^c_5(0,0)|^2_{\mathcal{ A}_5}
	\end{equation}
	for all $(t,j)\in \dom x^c_5$.	
\end{thm}
\begin{rem}
	 The proof of Theorem \ref{theo:theo_5} can be conducted using similar steps as in the proof of Theorem \ref{theo:theo_4}, which is omitted here. In view of (\ref{eqn:CRE}),  (\ref{eqn:gamma}), (\ref{eqn:set Q}) and (\ref{eqn:observer_X5})-(\ref{eqn:dx3}), one obtains the following hybrid closed-loop system:
	 \begin{equation}
	 \mathcal{H}^c_5:  \begin{cases}
	 \dot{x}^c_5 ~~= F_5(x^c_5),   & x^c_5\in \mathcal{F}^c_5 \\
	 {x^c_5}^+  = G_5(x^c_5),    &   x^c_5\in \mathcal{J}^c_5
	 \end{cases} \label{eqn:closed-loop5}
	 \end{equation}
	 where the flow and jump sets are defined as: $\mathcal{F}^c_5:=\{ x^c_5=(\hat{X}, \tilde{R}, {\mathsf{x}},t, \hat{b}_\omega ,\tilde{b}_\omega,\hat{b}_a,\tilde{b}_a) \in \mathcal{S}^c_5: \hat{X}\in \mathcal{F}_o \}$ and  $\mathcal{J}^c_5 :=\{x^c_5=(\hat{X}, \tilde{R}, {\mathsf{x}},t, \hat{b}_\omega ,\tilde{b}_\omega,\hat{b}_a,\tilde{b}_a) \in \mathcal{S}^c_5: \hat{X} \in \mathcal{J}_o\}$ with $\mathcal{F}_o$ and $\mathcal{J}_o$ given in (\ref{eqn:F_map}) and (\ref{eqn:J_map}), respectively. The flow and jump maps are given by
	 \begin{align*}
	 F_5(x^c_5)&= \begin{pmatrix}
	 f(\hat{X},\omega_y-\hat{b}_\omega,a_y-\hat{b}_a) - \Delta\hat{X}\\
	  \tilde{R} ((\hat{R}\tilde{b}_\omega)^\times - k_R\mathbb{P}_a(M\tilde{R}))\\
	 A(t)  {\mathsf{x}} - L C  {\mathsf{x}} + \nu \\
	 1\\
	 -k_\omega   \hat{R}\T \psi(M\tilde{R}) \\
	 -k_\omega   \hat{R}\T \psi(M\tilde{R}) \\
	 -k_c  \hat{R}\T  K_a  \tilde{R}\T \tilde{p}_e \\
	 -k_c  \hat{R}\T  K_a  \tilde{R}\T \tilde{p}_e
	 \end{pmatrix}  \\
	 G_5(x^c_5) &= \begin{pmatrix}
	 (X_q^{-1} \hat{X})\T,
	 (\tilde{R}R_q)\T,
	  {\mathsf{x}}\T ,
	 t ,
	 \hat{b}_\omega\T,
	 \tilde{b}_\omega\T,
	 \hat{b}_a\T ,
	 \tilde{b}_a\T
	 \end{pmatrix}\T.
	 \end{align*}
	 Note that the sets $\mathcal{F}^c_5$ and $\mathcal{J}^c_5$ are closed, and $\mathcal{F}^c_5 \cup \mathcal{J}^c_5 = \mathcal{S}^c_5$. Note also that the closed-loop system (\ref{eqn:closed-loop5}) satisfies the hybrid basic conditions of \cite{goebel2009hybrid} and is autonomous by taking $\omega_y$, $a_y$ and $L$ as functions of time $t$.
\end{rem}

\section{Simulation Results}\label{sec:simulation}

In this section, simulation results are presented to illustrate the performance of the proposed hybrid observers. We make use of the HyEQ Toolbox in Matlab \cite{sanfelice2013toolbox}. We refer to the continuous inertial navigation observer (\ie, observer $\mathcal{H}^o_3$ without jumps) as `CINO', the fixed-gain hybrid inertial navigation observer $\mathcal{H}^o_3$ as `HINO',  the CRE-based variable-gain hybrid inertial navigation observer $\mathcal{H}^o_4$ as `HINO-CRE', and the CRE-based variable-gain hybrid inertial navigation observer $\mathcal{H}^o_5$ as `HINO-CRE2'.

We consider an autonomous vehicle moving on a 10-meter diameter circle at 10-meter height, with the trajectory: $p(t)=10[\cos(0.8t),\sin(0.8t), 1]$. Consider the initial rotation as $R(0)=I_3$ and the angular velocity as $\omega(t) = [\sin(0.3\pi),0,0.1]\T$. Moreover, 6 landmarks are randomly selected such that Assumption \ref{assum:1} holds. We consider the same initial conditions for each observer as: $\hat{R}(0) = \mathcal{R}_a(0.99\pi,u), u\in \mathcal{E}(M), \hat{v}(0),\hat{p}(0)=\hat{b}_\omega = 0_{3\times 1}$.
We consider the gain parameters $k_i=1/6, i=1,2,\cdots,6, k_R=1, k_v=k_p=3, k_w =1$.  For the hybrid design, we choose $\mathbb{U}=\mathcal{E}(M)$, $\theta=0.8\pi$, and $\delta = 0.3 (1-\cos\theta)\Delta_M^*$ with $\Delta_M^*$ designed as per Lemma \ref{lemma:Delta_M}. Two sets of simulation are presented: the first one considers biased angular velocity measurements and the second one considers biased angular velocity and linear acceleration. The matrix parameters for the first CRE are chosen as $P(0) = 0.5I_5, V(t) = I_6, Q(t) = 10I_3$, and for the second CRE are chosen as $P(0) = I_9, V(t) = 0.05I_9, Q(t) = 10I_3$.

\begin{figure}[thpb]
	\centering
	\includegraphics[width=0.9\linewidth]{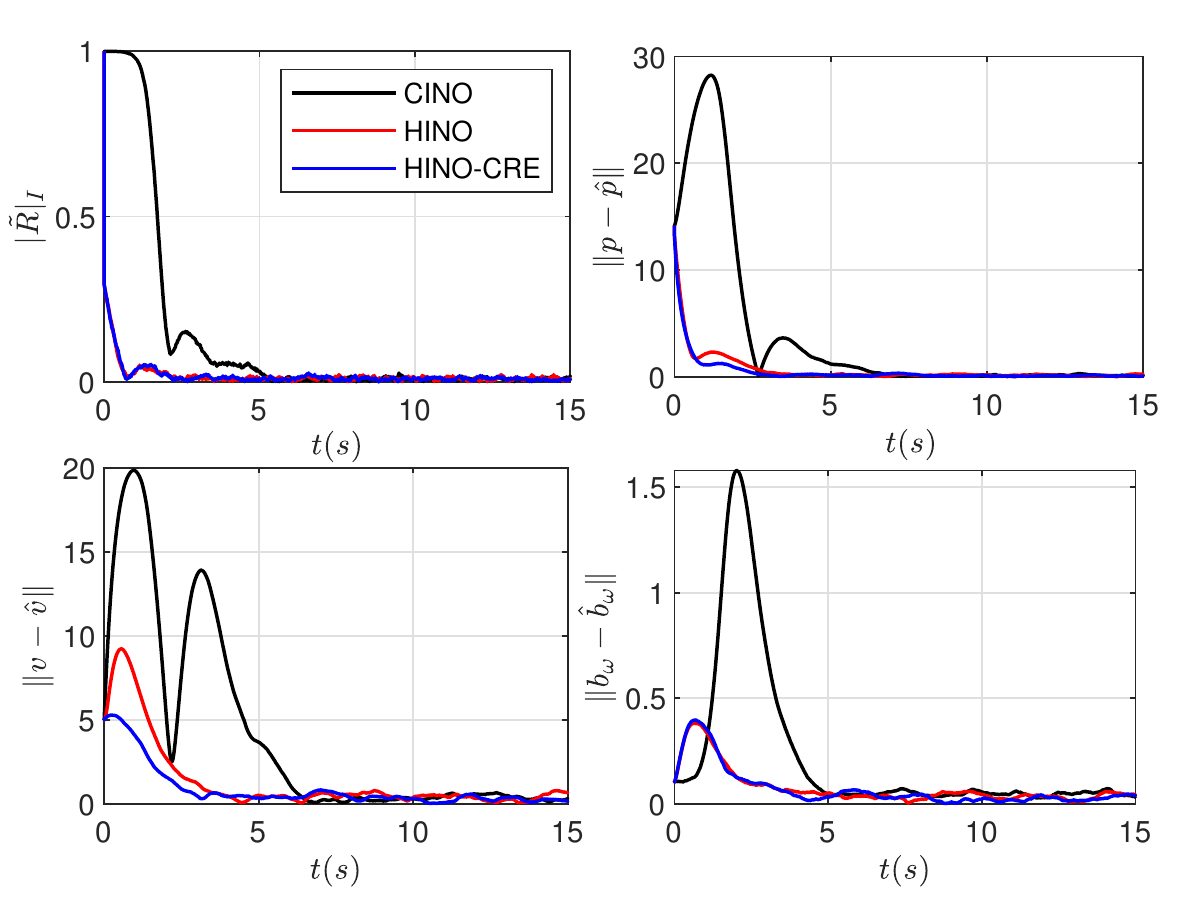}
	\caption{Simulation results with biased angular velocity $b_\omega=[-0.1~ 0.02~ 0.02]\T$ and additive white Gaussian noise of 0.4 variance in the measurements of $\omega$ and $a$, and 0.1 variance in the landmark position measurements. }
	\label{fig:simulation}
\end{figure}

\begin{figure}[thpb]
	\centering
	\includegraphics[width=0.95\linewidth]{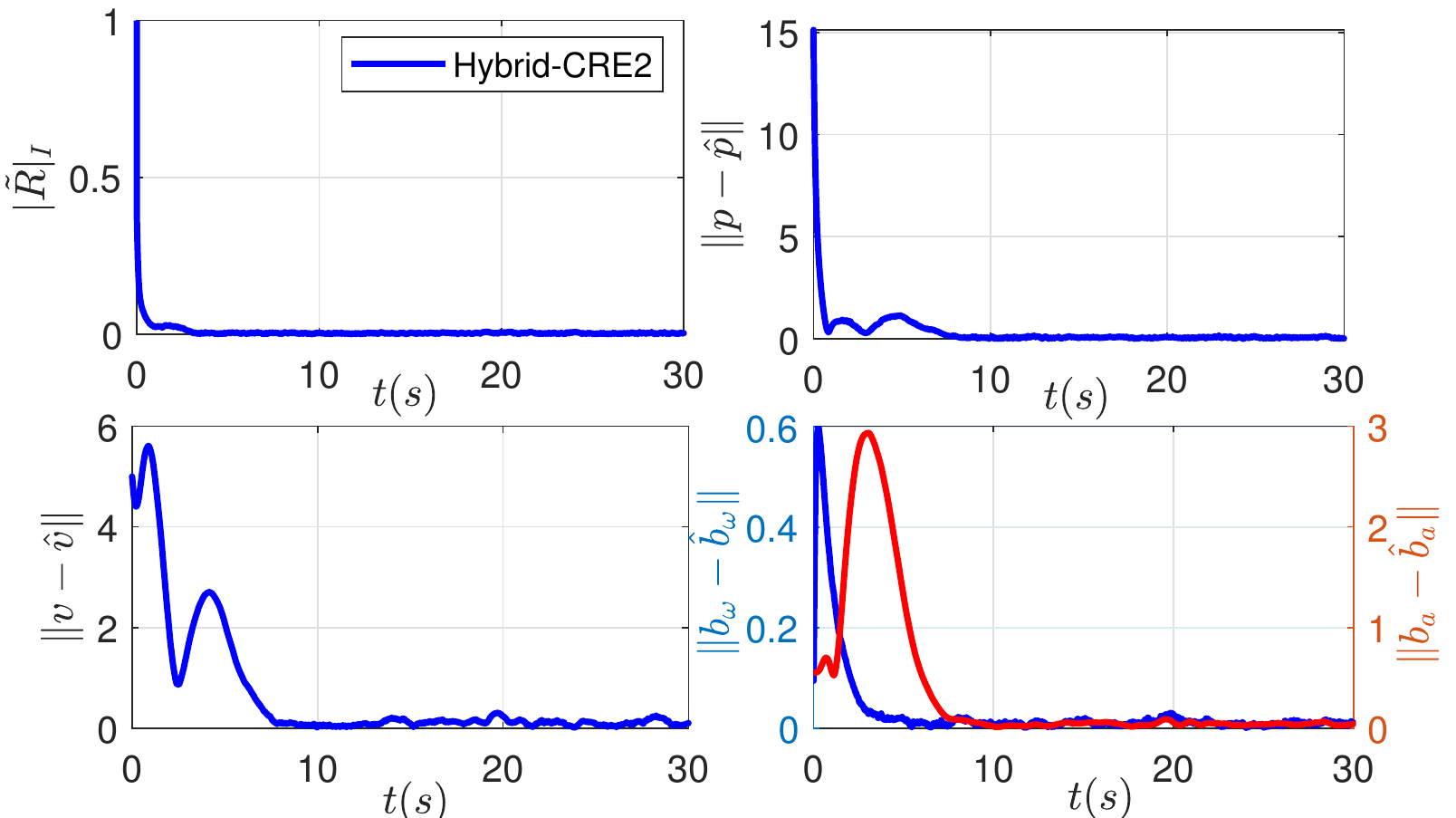}
	\caption{Simulation results with biased angular velocity $b_\omega=[-0.1~ 0.02~ 0.02]\T$, linear acceleration $b_a = [-0.01~	0.55 ~0.07]\T$, and additive white Gaussian noise of 0.1 variance in the angular velocity, linear acceleration and landmark position measurements. }
	\label{fig:simulation1}
\end{figure}

Simulation results are shown in Fig. \ref{fig:simulation} and Fig. \ref{fig:simulation1}. As one can see, the proposed hybrid observers exhibit fast convergence when the initial conditions are large.
Simulation results also illustrate the good performances of the proposed hybrid observers in the presence of biases in the angular velocity and linear acceleration.

\section{Experimental results}\label{sec:experimental}
To further validate the performance of our proposed hybrid observers, we applied our algorithms to real data from the EuRoc dataset \cite{Burri25012016}, where the trajectories are generated by a real flight of a quadrotor. This dataset includes a set of  images obtained from stereo cameras, IMU measurements and ground truth. The sampling rate of the IMU measurements from ADIS16448 is 200Hz and the sampling rate of the stereo images from MT9V034 is 20Hz. The ground truth of the states are obtained by a nonlinear least-squares batch solution using the Vicon pose and IMU measurements. More details about the EuRoC dataset can be found in \cite{Burri25012016}.

\subsection{Experimental setting}

\begin{figure}[thpb]
	\centering
	\includegraphics[width=0.9\linewidth]{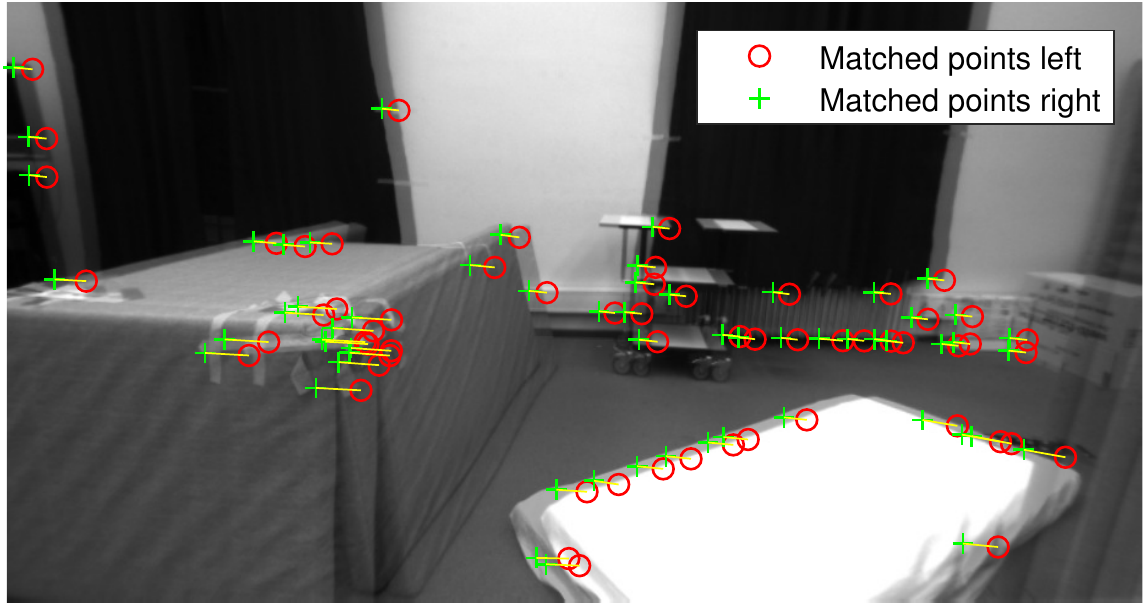}
	\caption{Example of features detection and tracking in the left and right images from a stereo camera using the MATLAB Computer Vision System Toolbox. Pictures come from the EuRoc dataset\cite{Burri25012016}.}
	\label{fig:features}
\end{figure}

The images are undistorted with the camera parameters calibrated using Stereo Camera Calibrator App in MATLAB.
The features are tracked via the Kanade-Lucas-Tomasi (KLT) tracker using minimum eigenvalue feature detection \cite{shi1994good}, which are shown in Fig. \ref{fig:features}. Since no physical landmarks are available in the EuRoc dataset, a set of `virtual' landmarks are generated from the stereo images and the ground truth pose at the beginning. More precisely, the coordinate of the $i$-th landmark expressed in the inertial frame is calculated as $p_i= R_{G}y_i + p_{G}$, where $R_G, p_G$ are the ground truth rotation and position of the vehicle, $y_i$ denotes the current three-dimensional position of the $i$-th point-feature generated from the current stereo images. For the sake of efficiency, we limit the maximum number of detected and tracked point-features to a certain number (60 in our experiments). It is quite unrealistic to track the same set of point-features through a long time image sequence. Hence, when the number of visible point-features is less than a certain threshold (6 in our experiments), a new set of point-features is generated using current stereo images and ground truth again.
The three-dimensional coordinates of the point-features from stereo images expressed in the camera frame (cam0) are transformed to the frame attached to the vehicle using the calibration matrix provided in the dataset. To remove matched point-feature outliers the technique proposed in \cite{hua2018attitude} has been used by choosing the thresholds $S=30, D=6$.

\subsection{Realtime implementation}
\begin{algorithm}\label{algo:1}
	\caption{Continuous-discrete algorithm of  HINO-CRE2}
	\begin{algorithmic}[1]
		\renewcommand{\algorithmicrequire}{\textbf{Initialization:}}
		\renewcommand{\algorithmicensure}{\textbf{Output:}}
		\REQUIRE  $\hat{X}(t_0)\in SE_2(3),\hat{b}_\omega(t_0)\in \mathbb{R}^3 ,\hat{b}_a(t_0)\in \mathbb{R}^3$, $P(t_0)\in \mathbb{R}^{9\times 9}>0$. 	
		\ENSURE  $\hat{X}(t), \hat{b}_{\omega}(t), \hat{b}_{a}(t)$ for all $t\geq t_0$
		\FOR {$1\leq k $} 	
		\WHILE {$ t_{k-1}\leq t\leq t_{k}$}		
		\STATE Integrate the following equations:
		\begin{flalign*}
		\begin{cases}
		\dot{\hat{X}} &= f(\hat{X},\omega_y-\hat{b}_\omega,a_y-\hat{b}_a)\\
		\dot{\hat{b}}_\omega &=0_{3\times 1}\\
		\dot{\hat{b}}_a &=0_{3\times 1}\\
		\dot{P} &= A(t) P + PA(t)\T + V(t)
		\end{cases}&&
		\end{flalign*}
		\ENDWHILE
		\STATE
		Set $\hat{X}_{k|k-1}=\hat{X}(t_{k}),~
		\hat{b}_{\omega,k|k-1}=\hat{b}_{\omega}(t_{k}),~ \hat{b}_{a,k|k-1}=\hat{b}_{a}(t_{k})$ and $
		P_{k|k-1}=P(t_{k})$ 	
		\STATE 	Compute the gain matrices
		\begin{flalign*}
		\begin{cases}
		L_k &= P_{k|k-1}C\T(CP_{k|k-1}C\T + Q^{-1}(t))^{-1}  \\
		K_p &=\frac{1}{k_c}  \hat{R}_{k|k-1} L_{1,k} \hat{R}_{k|k-1}\T\\
		K_v &= \frac{1}{k_c}  \hat{R}_{k|k-1} L_{2,k} \hat{R}_{k|k-1}\T\\
		K_a &= \frac{1}{k_c}  \hat{R}_{k|k-1} L_{3,k} \hat{R}_{k|k-1}\T
		\end{cases} &&
		\end{flalign*}
		from $L_k = [L_{1,k}\T, L_{2,k}\T, L_{3,k}\T]\T$ and $\hat{X}_{k|k-1} = \mathcal{T}(\hat{R}_{k|k-1},\hat{v}_{k|k-1},\hat{p}_{k|k-1})$
		\STATE  Compute the innovation terms   $\Delta_k$ in \eqref{eqn:innovation_term6} with $K_p$ and $K_v$,  $\Delta_{R,k}$ in (\ref{eqn:Delta_R}) and $\Delta_{p,k}$ in (\ref{eqn:Delta_p})	
		\STATE Update the state estimates as
		\begin{flalign*}
		\begin{cases}
		\hat{X}_{k|k} &=\exp(-\Delta_k)\hat{X}_{k|k-1} \\
		\hat{b}_{\omega,k|k} &= \hat{b}_{\omega,k|k-1} -  k_\omega \hat{R}_{k|k-1}\T \psi(\Delta_{R,k}) \\
		\hat{b}_{a,k|k} &= \hat{b}_{a,k|k-1} -    \hat{R}_{k|k-1}\T K_a \Delta_{p,k} \\
		P_{k|k} &= P_{k|k-1}-L_kCP_{k|k-1}
		\end{cases}&&
		\end{flalign*}
		\IF {($\mu_\mathbb{Q}(\hat{X}_{k|k},r,b_k) \geq \delta$)}
		\STATE Reset the state   $\hat{X}_{k|k} = X_q^{-1}\hat{X}_{k|k}, ~X_q\in \gamma(\hat{X}_{k|k})$
		\ENDIF
		\STATE   Set $\hat{X}(t_{k}) = \hat{X}_{k|k},~
		\hat{b}_{\omega}(t_{k})  = \hat{b}_{\omega, k|k}, ~\hat{b}_{a}(t_{k})  = \hat{b}_{a, k|k}$ and $
		P(t_{k})=P_{k|k}$
		\ENDFOR
	\end{algorithmic}
\end{algorithm}

In practice, the IMU measurements can be obtained at a high rate, while the landmark measurements are often obtained, for example with stereo cameras, at a much lower rate.  Taking into account this fact, we define a strictly increasing sequence $\{t_k\}_{k\in \mathbb{N}/\{0\}}$ as the time-instants when the landmark measurements are obtained. Inspired by the work of continuous-discrete Kalman filter and extended Kalman filter in \cite[page 194]{Lewis2007optimal} and \cite{kulikov2014accurate}, we implement our hybrid observer HINO-CRE2 as shown in Algorithm 1. The proposed algorithm has two parts: the states are continuously updated from IMU when no measurements of landmarks are received (\ie, $t\in(t_{k-1},t_{k})$); when the measurements arrive (\ie, $t=t_k$) the state variables are updated using the landmark measurements. This type of continuous-discrete observers for inertial navigation has also been considered in \cite{barrau2017invariant,hua2018attitude}. The CRE is continuously integrated from the time $t_{k-1}$ to the next time instant $t_k$ when the new landmark measurements arrive, and then a numerical discretization method is applied at the instant of time $t_k$. A first-order numerical discretization method is applied to the dynamics of the estimated angular velocity bias $\hat{b}_\omega$ and linear acceleration bias $\hat{b}_a$. However, an exponential map based discrete update of $\hat{X}$ has been considered, \ie, $\hat{X}_{k|k} = \exp(-\Delta_k) \hat{X}_{k|k-1}$, which guarantees that $\hat{X}_{k|k} \in SE_2(3)$.
Note that the estimated state $\hat{X}$ is reset once the condition $\mu_\mathbb{Q}(\hat{X},r,b) \geq \delta$ (\ie, $\hat{X}\in \mathcal{J}_o$) is satisfied. The Algorithm 1 can be easily modified for the other observers proposed in this paper which are omitted here.

\subsection{Results}

\begin{figure*}[thpb]
	\centering
\begin{minipage}{0.95\linewidth}
		\includegraphics[width=0.34\linewidth]{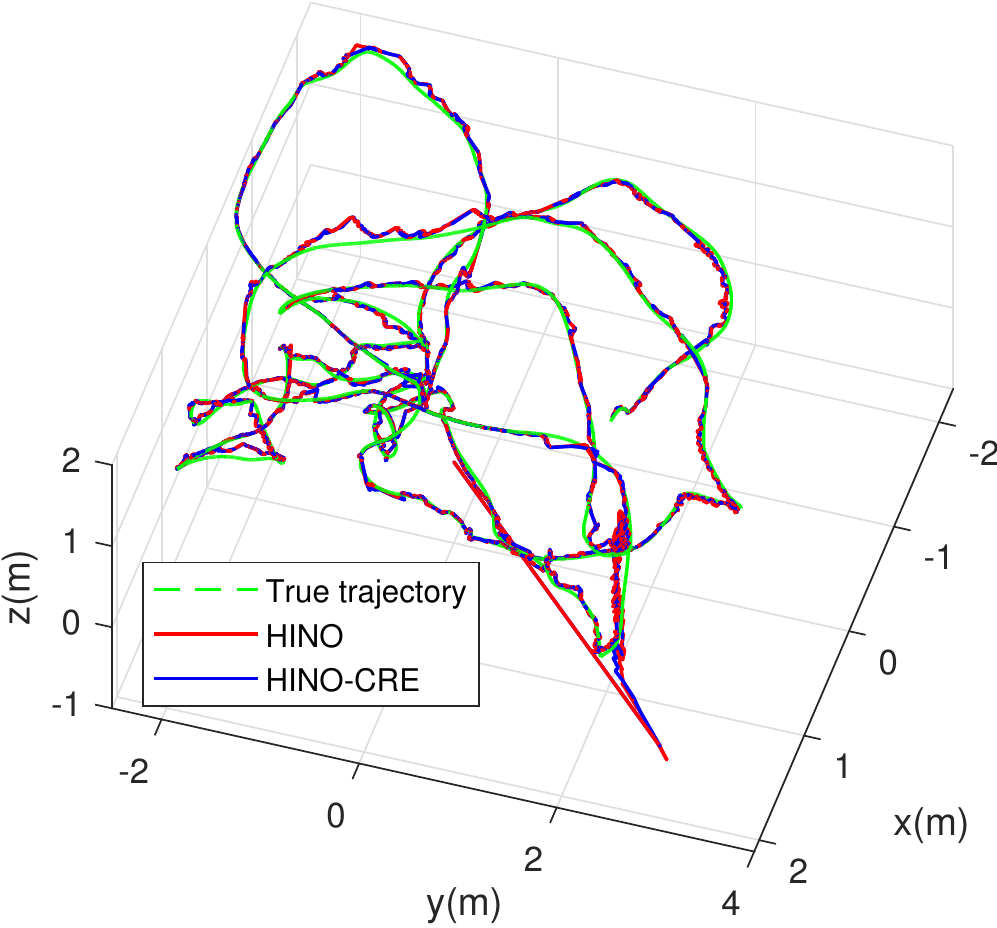}
	\includegraphics[width=0.6\linewidth]{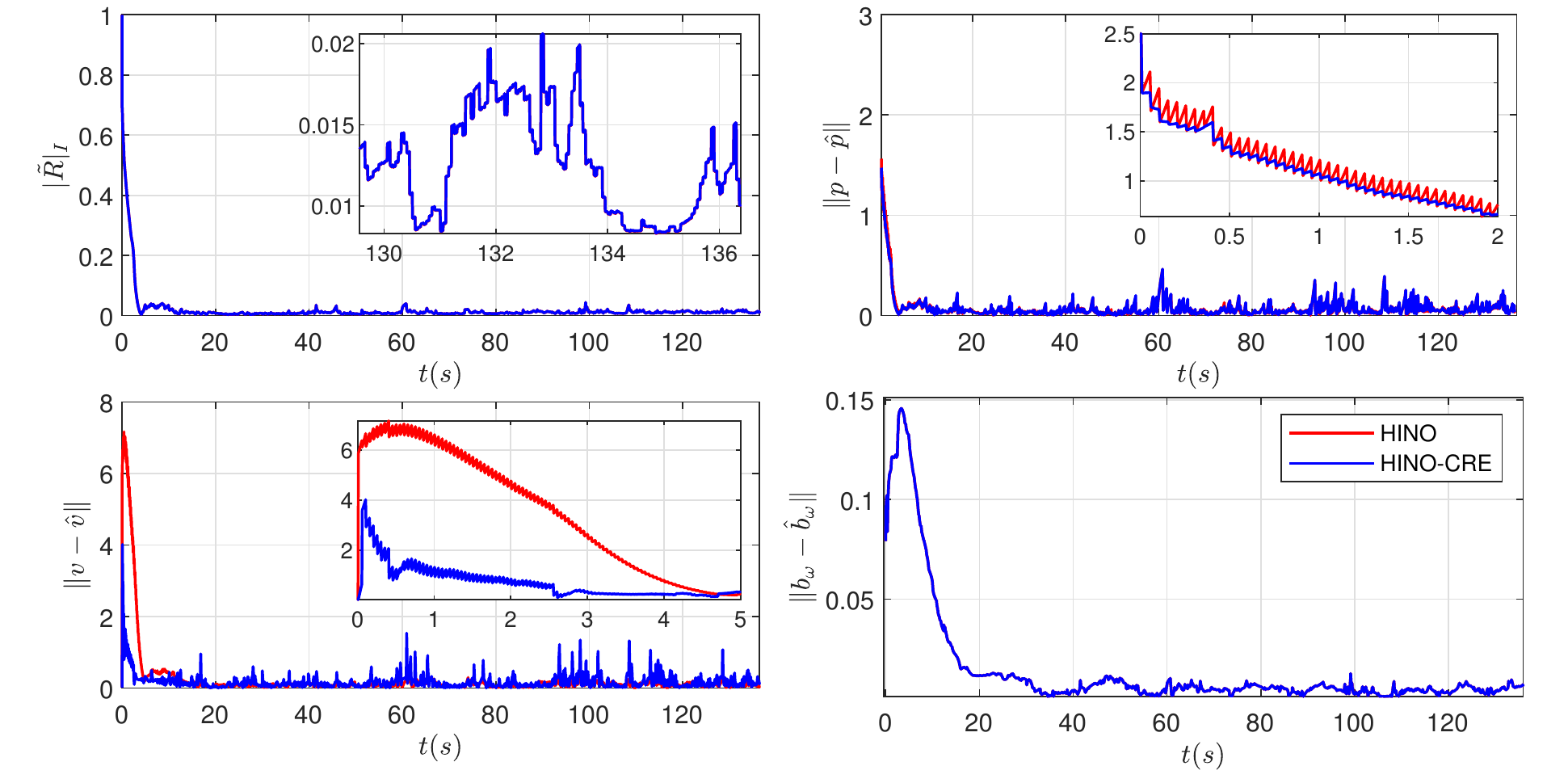}
\end{minipage} \caption*{(a) Experimental results using biased gyro and unbiased accelerometer measurements}
\begin{minipage}{0.95\linewidth}
	\includegraphics[width=0.34\linewidth]{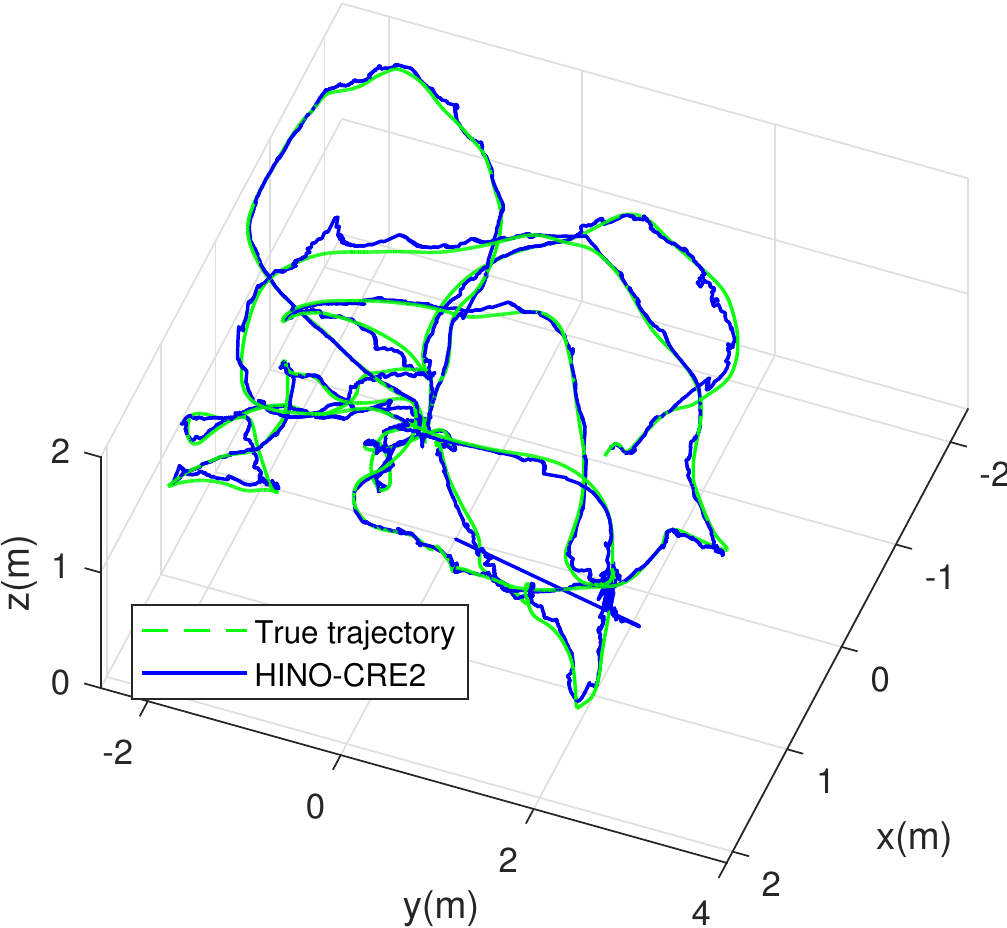}
	\includegraphics[width=0.63\linewidth]{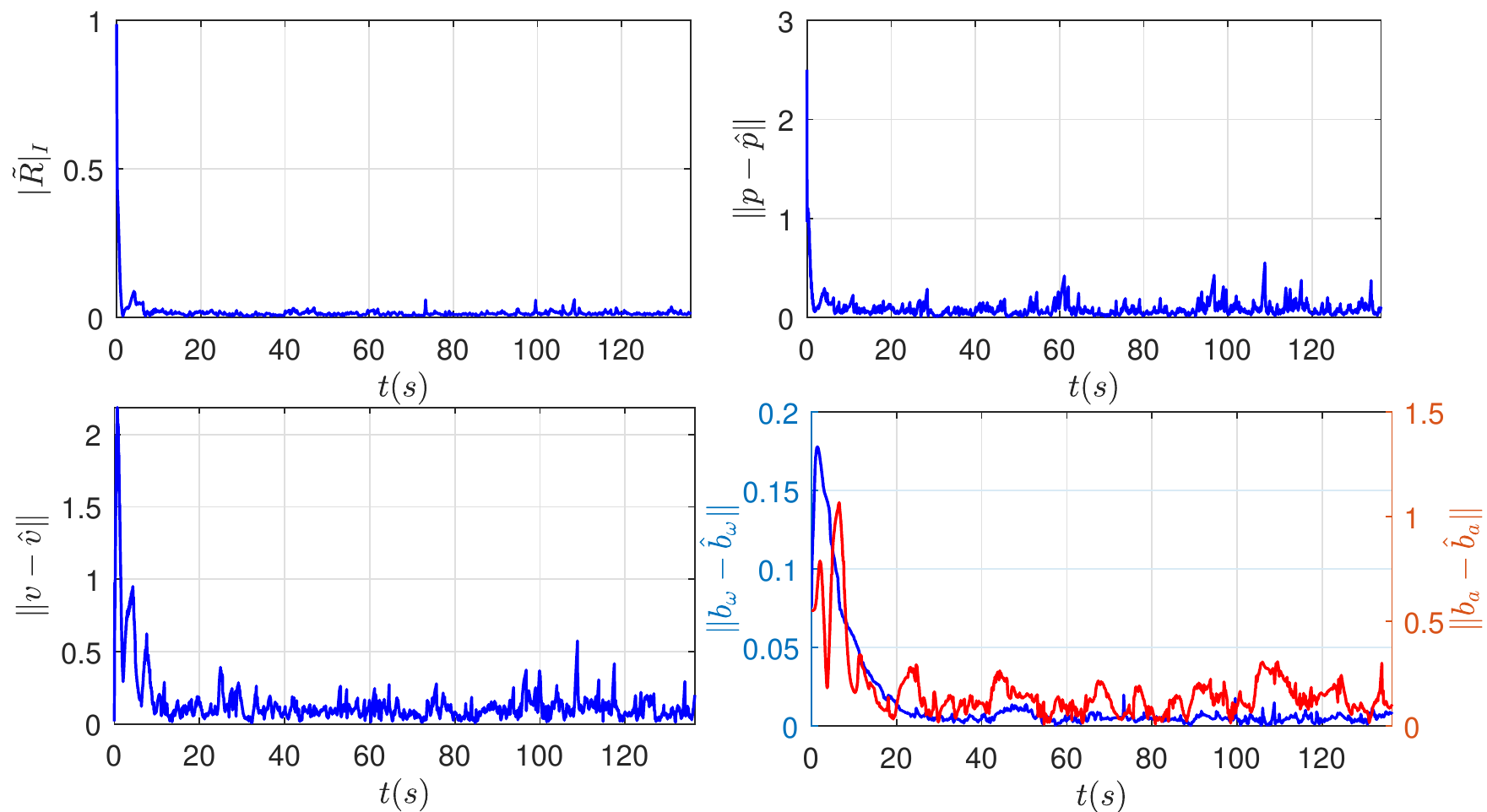}
\end{minipage}\caption*{(b) Experimental results using biased gyro and accelerometer measurements}
	\caption{The experimental results use the data of Vicon Room 1 01 from EuRoc dataset\cite{Burri25012016} with large initial conditions: $\hat{R}(0)=\exp(0.99\pi e_3^\times)R_G$ and $ \hat{p}(0)=\hat{v}(0)=\hat{b}_\omega(0)=\hat{b}_a(0) = 0_{3\times 1}$.
		The true and estimated trajectories are shown in the left plot. The estimation errors of rotation, position, velocity and IMU bias are shown in the right plot.
	}
	\label{fig:experiment}
\end{figure*}

Two sets of experiments have been presented in this paper, and
the gain parameters are carefully tuned with a trade-off between the convergence rate and the noise at steady state. Note that higher gains result in faster convergence but amplify noise at steady-state. As one can see, in the first plot of Fig. \ref{fig:experiment}, the estimates provided by both observers HINO and HINO-CRE, using the biased gyro and unbiased accelerometer measurements from IMU and stereo vision, converge, after a few seconds, to the vicinity of the ground truth. In the second plot of Fig.\ref{fig:experiment}, the estimates (including accelerometer-bias) provided by observer HINO-CRE2, using the biased IMU measurements and stereo vision, also converge to the vicinity of the ground truth after a few seconds. Note that the ground truth pose is used to validate the performance of the proposed algorithms and also to generate the virtual landmarks in the experiments due to the lack of physical landmarks.

\section{Conclusions}

Nonlinear geometric hybrid observers for inertial navigation systems, with global exponential stability guarantees, have been proposed. The observers are designed on the matrix Lie group $SE_2(3)$ using IMU and landmark position measurements, relying on a resetting mechanism designed to avoid the undesired equilibrium points in the flows and to ensure a decrease of the Lypunov function after each jump. Both non-biased and biased IMU measurements have been considered. Variable-gain versions of these observers, relying on a CRE, have also been proposed to efficiently handle measurements noise. Simulation and experimental results, illustrating the performance of the proposed hybrid observers, have been provided. For future work, it will be interesting to investigate the practical scenario of multi-rate intermittent measurements.


\appendix

\subsection{Proof of Proposition \ref{pro:undesired_eq}} \label{sec:undesired_eq}
For each $x^c_1\in SE_2(3)\times \Psi \times \mathbb{R}^+$, let us rewrite $\tilde{R}=\mathcal{R}_a(\pi,v)$ with $v\in \mathcal{E}(M)$, and $ R_q = \mathcal{R}_a(\theta,u_q)$ with $\theta\in(0,\pi]$ and $u_q\in \mathbb{U}$. In view of (\ref{eqn:Phi_R}) and (\ref{eqn: definition_mu_R}), one can show that
\begin{align*}
\mu_{\mathbb{Q}}(\hat{R},r,b)
&= \Upsilon(\hat{X},r,b) - \min_{X_q=\mathcal{T}(R_q,p_q,v_q)\in \mathbb{Q}} \Upsilon(X_q\T \hat{X},r,b)   \\
&= \tr((I_3-\tilde{R})M) - \min_{R_q\in \mathcal{R}_a(\theta,\mathbb{U})}\tr((I_3-\tilde{R}R_q)M)\\
&=\max_{R_q\in \mathcal{R}_a(\theta,\mathbb{U})} \tr(\tilde{R}(I_3-R_q)M)\\
& =   (1-\cos\theta)\max_{u_q\in \mathbb{U}} \Delta(u_q,v) \\
&   \geq (1-\cos\theta) \max_{u_q\in \mathcal{E}(M)} \Delta(u_q,v)
\end{align*}
where we made use of the definition (\ref{eqn:definition_Delta_uv}) and the fact that $\max_{u_q\in \mathbb{U}} \Delta(u_q,v)\geq \max_{u_q\in \mathcal{E}(M)} \Delta(u_q,v)$ for any $v\in \mathbb{R}^3$. From the definition of $\Delta_M^*$ given in (\ref{eqn:Delta_M_star}) such that  for any $x^c_1\in SE_2(3)\times \Psi \times \mathbb{R}^+ $, one has
\begin{align*}
\mu_{\mathbb{Q}}(\hat{R},r,b)
&  \geq (1-\cos\theta) \min_{v\in \mathcal{E}(M)} \max_{u_q\in \mathcal{E}(M)} \Delta(u,v)  \\
&\geq (1-\cos\theta) \Delta^*_M > \delta
\end{align*}
which gives $SE_2(3)\times \Psi \times \mathbb{R}^+ \subseteq \mathcal{J}^c_1$ from (\ref{eqn:J_map}) and (\ref{eqn:closed-loop}). This completes the proof.

\subsection{Proof of Theorem \ref{theo:theo_1}}\label{sec:theo_1}

From Proposition \ref{pro:pro_1} and \ref{pro:undesired_eq}, one can show that all the undesired equilibrium points of the flows of $\mathcal{H}^c_1$ lie in the jump set $\mathcal{J}^c_1$.
Consider the following real-valued function $\mathcal{L}_R:SO(3)\to \mathbb{R}^+$:
\begin{equation}
\mathcal{L}_R(\tilde{R}) = \tr((I-\tilde{R})M)   \label{eqn:L_R}.
\end{equation}
Let $\bar{M}: = \frac{1}{2}(\tr(M)I_3-M), \underline{M} := \tr(\bar{M}^2)I - 2\bar{M}^2$ and $\underline{\bar{M}}: = \frac{1}{2}(\tr(\underline{M})I_3-\underline{M})$. Applying the results in  \cite[Lemma 2]{berkane2017hybrid}, one obtains
\begin{align}
 4\lambda_{m}^{\bar{M}} |\tilde{R}|_I^2 &\leq \mathcal{L}_R \leq  4\lambda_{M}^{\bar{M}}  |\tilde{R}|_I^2 \label{eqn:L_R_bound}, \\
 \dot{\mathcal{L}}_R &\leq -\lambda_R |\tilde{R}|_I^2 \qquad  \qquad    \forall x^c_1 \in \mathcal{F}^c_1 \label{eqn:dot_L_R}
\end{align}
where  $\lambda_R: =  4k_R\varrho_M \lambda_{m}^{\underline{\bar{M}}} $, and $\varrho_M: = \min_{ x^c_1 \in \mathcal{F}^c_1} \varrho(M,\tilde{R}) $ with $\varrho(M,\tilde{R}) :=(1-|\tilde{R}|_I^2\cos^2(u,\bar{M}u)) $ and $u\in \mathbb{S}^2$ denoting the axis of rotation $\tilde{R}$. Moreover, one can verify that for all $x_1^c\in \mathcal{F}^c_1$ one has $\varrho(M,\tilde{R})>0$ which implies $\lambda_R>0$ in the flows.

On the other hand, consider the following real-valued function $\mathcal{L}_p: \mathbb{R}^3\times \mathbb{R}^3 \to \mathbb{R}^+$ as
\begin{equation}
\mathcal{L}_p(\tilde{p}_e,\tilde{v}) = \frac{1}{2} \|\tilde{p}_e\|^2 + \frac{1}{2 k_c k_v} \|\tilde{v}\|^2 - \mu \tilde{p}_e\T \tilde{v}  \label{eqn:Lp}
\end{equation}
with some $\mu>0$.
Let $e_2 := [\|\tilde{p}_e\|~ \|\tilde{v}\|]\T$. One verifies that
\begin{align}
e_2\T \underbrace{\begin{bmatrix}
	\frac{1}{2} & -\frac{\mu}{2} \\
	-\frac{\mu}{2} & \frac{1}{2  k_c k_v}
	\end{bmatrix}}_{P_1}e_2  \leq \mathcal{L}_p \leq e_2\T  \underbrace{\begin{bmatrix}
	\frac{1}{2} & \frac{\mu}{2} \\
	\frac{\mu}{2} & \frac{1}{2  k_c k_v}
	\end{bmatrix}}_{P_2} e_2  \label{eqn:L_p_bound}.
\end{align}
The time-derivative of $\mathcal{L}_p$ along the flows of (\ref{eqn:closed-loop}) is given by
\begin{align*}
\dot{\mathcal{L}}_p & =  \tilde{p}_e\T(-k_c k_p \tilde{p}_e + \tilde{v}) + \frac{1}{ k_c k_v} \tilde{v}\T (-k_c k_v \tilde{p}_e  + (I-\tilde{R})\mathsf{g}) \\
&\quad - \mu(-k_c k_p \tilde{p}_e+\tilde{v})\T \tilde{v} -\mu \tilde{p}_e\T (-  k_c k_v \tilde{p}_e + (I_3-\tilde{R})\mathsf{g})\\
& =  -k_c k_p \tilde{p}_e\T \tilde{p}_e + \mu  k_c k_v \tilde{p}_e\T \tilde{p}_e - \mu \tilde{v}\T \tilde{v} + \mu k_c k_p \tilde{p}_e\T \tilde{v}  \\
& \quad+  \frac{1}{  k_c k_v} \tilde{v}\T   (I-\tilde{R})\mathsf{g}-  \mu \tilde{p}_e\T   (I_3-\tilde{R})\mathsf{g}\\
& \leq    -(k_p-\mu k_v ) k_c \|\tilde{p}_e\|^2  - \mu \|\tilde{v}\| + \mu k_c k_p \|\tilde{p}_e\| \|\tilde{v}\|  \\
& \quad+  \frac{1}{  k_c k_v} \|\g\|\|\tilde{v}\| \|I-\tilde{R}\|_F + \mu \|\g\| \|\tilde{p}_e\| \|I_3-\tilde{R}\|_F .
\end{align*}
Let $c_1: = \max\{ \frac{\|\g\|}{  k_c k_v} , \mu \|\g\|\}$, one can further deduce that
\begin{align}
\dot{\mathcal{L}}_p & \leq  - e_2\T \underbrace{\begin{bmatrix}
	(k_p-\mu k_v ) k_c & \frac{\mu k_c k_p}{2} \\
	\frac{\mu k_c k_p}{2} & \mu
	\end{bmatrix} }_{P_3}  e_2 + 4c_1|\tilde{R}|_I \|e_2\|  \label{eqn:dot_L_p}
\end{align}
where we made use of the facts: $ \|I_3-\tilde{R}\|_F = 2\sqrt{2}|\tilde{R}|_I$ and $(\|\tilde{v}\| + \|\tilde{p}_e\|) \leq \sqrt{2(\|\tilde{v}\|^2 + \|\tilde{p}_e\|)^2} =\sqrt{2}\|e_2\|$.
To guarantee that the matrices $P_1,P_2$ and $P_3$ are positive definite, it is sufficient to pick $\mu$ as
$$
0< \mu < \min\left\{\frac{1}{\sqrt{k_c k_v}}, \frac{4k_p}{4k_v + k_ck_p^2} \right\}.
$$

To show exponential stability, let us consider the following Lyapunov function candidate:
\begin{equation}
\mathcal{L}(x^c_1): = \mathcal{L}_R(\tilde{R}) + \varepsilon \mathcal{L}_p(\tilde{p}_e,\tilde{v})  \label{eqn:Lyp}
\end{equation}
with some $0<\varepsilon$.
Let $|x^c_1|_{\mathcal{A}_1}\geq 0$ denote the distance to the set $\mathcal{A}_1$ such that
$
|x^c_1|_{\mathcal{A}_1}^2  : = \textstyle \inf_{y=(\bar{X},I_3,0,0,\bar{t})\in \mathcal{A}_1}  (\|\bar{X}-\hat{X}\|_F^2 + |\tilde{R}|_I^2
 + \|\tilde{p}_e\|^2 + \|\tilde{v}\|^2+\|\bar{t}-t\|^2 )
=   |\tilde{R}|_I^2 + \|\tilde{p}_e\|^2 + \|\tilde{v}\|^2 =  |\tilde{R}|_I^2 + \|e_2\|^2 .
$
From (\ref{eqn:L_R_bound}) and (\ref{eqn:L_p_bound}), one has
\begin{align}
\underline{\alpha} |x^c_1|_{\mathcal{A}_1}^2 \leq \mathcal{L}(x^c_1) \leq \bar{\alpha} |x^c_1|_{\mathcal{A}_1}^2     \label{eqn:L_x}
\end{align}
where $\underline{\alpha}: = \min\{4\lambda_{m}^{\bar{M}},\varepsilon \lambda_{m}^{P_1}\}, \bar{\alpha}: =\max\{4\lambda_{M}^{\bar{M}},\varepsilon\lambda_{M}^{P_2}\}$. From (\ref{eqn:dot_L_R}) and (\ref{eqn:dot_L_p}), one has
\begin{align}
\dot{\mathcal{L}}(x^c_1) & \leq -\lambda_R |\tilde{R}|_I^2 - \varepsilon\lambda_{m}^{P_3}\|e_2\|^2 + 4\varepsilon c_1|\tilde{R}|_I \|e_2\|  \nonumber \\
& = - \begin{bmatrix}
|\tilde{R}|_I & \|e_2\|
\end{bmatrix} \underbrace{\begin{bmatrix}
	\lambda_R & -2\varepsilon c_1\\
	-2\varepsilon c_1 & \varepsilon\lambda_{m}^{P_3}
	\end{bmatrix}}_{P_4}  \begin{bmatrix}
|\tilde{R}|_I \\ \|e_2\|
\end{bmatrix} \nonumber \\
&\leq -\lambda_F  \mathcal{L}(x^c_1)  \label{eqn:dL_x}
\end{align}
where $P_4$ is positive definite by choosing $\varepsilon<{\lambda_R\lambda_{m}^{P_3}}/{(4c_1^2)}$, and $\lambda_F:= \lambda_{m}^{P_4}/\bar{\alpha}$ with $\bar{\alpha}$ given in (\ref{eqn:L_x}).  In view of the jumps of  (\ref{eqn: definition_mu_R})-(\ref{eqn:J_map}), (\ref{eqn:closed-loop}) and (\ref{eqn:L_R}), one shows	
\begin{align}
& \mathcal{L}({x^c_1}^+) - \mathcal{L}(x^c_1) \nonumber \\
&= \mathcal{L}_R(\tilde{R}^+)  - \mathcal{L}_R(\tilde{R}) + \varepsilon \mathcal{L}_p(\tilde{p_e}^+,\tilde{v}^+)  - \varepsilon\mathcal{L}_p (\tilde{p_e},\tilde{v}) \nonumber \\
&= -\Upsilon(\hat{X},r,b) + \min_{ X_q \in \mathbb{Q}}   \Upsilon(X_q\T\hat{X},r,b)   \nonumber  \\
&= -\mu_\mathbb{Q}(\hat{X},r,b) \leq  -\delta  \label{eqn:L_x+}
\end{align}
where we made use of the facts: $\mathcal{L}_R=\Upsilon(\hat{X},r,b), \mathcal{L}_R^+=\min_{X_q\in \mathbb{Q}}   \Upsilon(X_q\T\hat{X},r,b)$ from (\ref{eqn:Phi_R})-(\ref{eqn: definition_mu_R}), and $\mathcal{L}_p^+ =\mathcal{L}_p$ from $\tilde{p}_e=\tilde{p}_e, \tilde{v}^+=\tilde{v}$.
Using the facts $\mathcal{L}_R^+  - \mathcal{L}_R \leq -\delta$ and (\ref{eqn:dot_L_R}),  one has
$
\mathcal{L}_R(\tilde{R}(t,j))
\leq \cdots  \leq \mathcal{L}_R(\tilde{R}(0,0))  - j \delta
$,
where $(t,j) \in  \dom x^c_1$. From (\ref{eqn:L_R_bound}), one obtains
$
j \leq J:=  \lceil  {4 \lambda_{M}^{\bar{M}} }/{\delta}  \rceil,
$
where $\lceil \cdot\rceil$ denotes the ceiling function. Hence, one can conclude that the number of jumps is finite.

Since the solution of $x^c_1$ is complete and the number of jumps is bounded, the hybrid time domain takes the form $\dom x^c_1 = \cup_{j=0}^{J-1}([t_j,t_{j+1}]\times \{j\})\cup ([t_J,+\infty) \times \{J\})$. In view of (\ref{eqn:dL_x})-(\ref{eqn:L_x+}), one obtains
$
\mathcal{L}(x^c_1(t,j)) \leq \exp(-\lambda_F (t-t_j)) \mathcal{L}(x^c_1(t_j,j))  \leq \exp(-\lambda_F t) \mathcal{L}(x^c_1(0,0))
 \leq \exp(\lambda_F J) \exp(-\lambda_F(t+j))   \mathcal{L}(x^c_1(0,0)).
$
Substituting (\ref{eqn:L_x}), one concludes that for each $(t,j)\in \dom x^c_1$,
\begin{equation*}
|x^c_1(t,j)|^2_{\mathcal{A}_1} \leq \kappa\exp\left(- \lambda_F(t+j)\right) |x^c_1(0,0)|^2_{\mathcal{ A}_1}
\end{equation*}
where $ \kappa:=\exp( \lambda_F J) \bar{\alpha}/ \underline{\alpha} $.  This completes the proof.

\subsection{Proof of Lemma \ref{lem:existenceP}}\label{sec:existenceP}
To show the uniform observability of the pair $(A(t),C)$, we need to explicitly calculate the state transition matrix $\Phi(t,\tau)$.
Consider a time-varying rotation matrix $\bar{R}(t)$ with $\bar{R}(0)\in SO(3)$ and $\dot{\bar{R}}(t) =  (-\omega(t))^\times \bar{R}(t) $. Note that $\bar{R}(t)$ does not have to be equal to $R(t)$. Inspired from \cite{hamel2018riccati}, let us introduce the matrices $T(t) = \text{blkdiag}([
\bar{R}(t),  \bar{R}(t)
])$, $ S(t) = \text{blkdiag}([
(-\omega(t))^\times, (-\omega(t))^\times
])$  and constant  matrix $\bar{A} = A(t) - S(t)$. Then, one verifies that $\dot{T}(t) =  S(t)T(t)$ and $T(t)\bar{A} = \bar{A}T(t)$.

Next, we are gonging to show that the state transition matrix $\Phi(t,\tau)$  can be expressed as
\begin{equation}
\Phi(t,\tau) =  T(t) \bar{\Phi}(t,\tau) T^{-1}(\tau)  \label{eqn:PhiA}
\end{equation}
with   $\bar{\Phi}(t,\tau) = \exp{(\bar{A}(t-\tau))}$ being the state transition matrix associated to $\bar{A}$. From (\ref{eqn:PhiA}), one can easily verify that $	\Phi(t,t)  = I_n, \Phi^{-1}(t,\tau) = \Phi(\tau,t) $ and $\Phi(t_3,t_2)\Phi(t_2,t_1)  = \Phi(t_3,t_1)$ for every $t_1,t_2,t_3\geq 0$. Applying the facts $\dot{T}(t) = S(t)T(t)$ and $T(t)\bar{A} = \bar{A}T(t)$, one can further show that
\begin{align*}
\frac{d}{dt}\Phi(t,\tau) & = \dot{T}(t) \bar{\Phi}(t,\tau) T^{-1}(\tau) + T(t) \frac{d}{dt}\bar{\Phi}(t,\tau) T^{-1}(\tau) \\
& = S(t)T(t) \bar{\Phi}(t,\tau) T^{-1}(\tau) + T(t) \bar{A}\bar{\Phi}(t,\tau) T^{-1}(\tau) \\
& = S(t)T(t) \bar{\Phi}(t,\tau) T^{-1}(\tau) + \bar{A}T(t) \bar{\Phi}(t,\tau) T^{-1}(\tau) \\
& = A(t)T(t) \bar{\Phi}(t,\tau) T^{-1}(\tau) \\
& = A(t)  {\Phi}(t,\tau).
\end{align*}
Therefore, one can conclude that $\Phi(t,\tau)$ is the state transition matrix  associated to $A(t)$. Using the facts $T^{-1}(\tau) = T(\tau)\T$, $T^{-1}(\tau) C \T = C \T \bar{R}(\tau)\T$ and $\bar{R}(\tau)\T \bar{R}(\tau)  = I_3$ one obtains
\begin{align}
&   \frac{1}{\delta}  \int_{t}^{t+\delta} T(t) \bar{\Phi}(\tau,t)\T  T^{-1}(\tau) C \T C  T(\tau) \bar{\Phi}(\tau,t) T^{-1}(t) d\tau \nonumber \\
&= T(t) \left( \frac{1}{\delta} \int_{t}^{t+\delta}  \bar{\Phi}(\tau,t)\T   C \T C   \bar{\Phi}(\tau,t)  d\tau \right) T^{-1}(t) \label{eqn:TPhiT}.
\end{align}	
Note that the pair $(\bar{A}, {C})$ is (Kalman) observable, \ie, $\text{rank}[ {C}, {C}\bar{A},\cdots, {C}\bar{A}^6] = 6$, and there exist positive constants $\bar{\delta},\bar{\mu}$ such that for all $t\geq 0$ one has
$
\bar{W}(t,t+\bar{\delta}) =
\frac{1}{\bar{\delta}} \int_{t}^{t+\bar{\delta}} \bar{\Phi} (\tau,t)\T  {C}\T  {C} \bar{\Phi}(\tau,t) d\tau \geq \bar{\mu} I_6 .
$
From (\ref{eqn:TPhiT}), choosing $\delta\geq \bar{\delta}$ and $0<\mu \leq \frac{\bar{\delta}}{\delta} \bar{\mu}$, it follows that
$$
W(t,t+\tau)
 \geq \frac{\bar{\delta}}{\delta} T (t)\bar{W}(t,t+\bar{\delta}) T^{-1}(t)
 \geq  \frac{\bar{\delta}}{\delta} I_6
 \geq \mu I_6
$$
for all $t\geq 0$, which implies that the pair $(A(t),C)$ is uniformly observable. This completes the proof.

\subsection{Proof of Theorem \ref{theo:theo_2}} \label{sec:theo_2}
The proof of Theorem \ref{theo:theo_2} is similar to the proof of Theorem \ref{theo:theo_1}.   In view of (\ref{eqn:CRE}), (\ref{eqn:set Q}),  (\ref{eqn:gamma}), and (\ref{eqn:observer_X2})-(\ref{eqn:dx}), one obtains the following hybrid closed-loop system:
\begin{equation}
\mathcal{H}^c_2:  \begin{cases}
\dot{x}^c_2 ~~= F_2(x^c_2),   & x^c_2\in \mathcal{F}^c_2 \\
{x^c_2}^+  = G_2(x^c_2),    &   x^c_2\in \mathcal{J}^c_2
\end{cases} \label{eqn:closed-loop2}
\end{equation}
where the flow and jump sets are defined as $\mathcal{F}^c_2:=\{ x^c_2=(\hat{X},\tilde{R},\mathsf{x},t) \in \mathcal{S}^c_2: \hat{X} \in \mathcal{F}_o,  \}$ and  $\mathcal{J}^c_2 :=\{ x^c_2=(\hat{X},\tilde{R},\mathsf{x},t) \in \mathcal{S}^c_2: \hat{X} \in \mathcal{J}_o \}$, and the flow and jump maps are given by
\begin{align*}
F_2(x^c_2) = \begin{pmatrix}
f(\hat{X},\omega,a) - \Delta\hat{X}\\
\tilde{R} (-k_R\mathbb{P}_a(M\tilde{R}))\\
A \mathsf{x} - LC \mathsf{x} + \nu \\
1
\end{pmatrix},
G_2(x^c_2) = \begin{pmatrix}
X_q^{-1} \hat{X} \\
\tilde{R}R_q\\
\mathsf{x}\\
t
\end{pmatrix}.
\end{align*}
Note that the sets $\mathcal{F}^c_2, \mathcal{J}^c_2$ are closed, and $\mathcal{F}^c_2 \cup \mathcal{J}^c_2 = \mathcal{S}^c_2$. Note also that the closed-loop system (\ref{eqn:closed-loop2}) satisfies the hybrid basic conditions of \cite{goebel2009hybrid} and is autonomous by taking $\omega$, $a$, $A $ and  $L$ as functions of $t$. By following similar steps as in the proof of Proposition \ref{pro:pro_1} and \ref{pro:undesired_eq}, one can show that all the undesired equilibrium points of the flows of $\mathcal{H}^c_2$ lie in the jump set $\mathcal{J}^c_2$, which is omitted here.

Consider the following Lyapunov function candidate:
\begin{equation}
\mathcal{L}(x^c_2): = \mathcal{L}_R(\tilde{R}) + \varepsilon \bar{\mathcal{L}}_p(\mathsf{x})  \label{eqn:Lyp2}
\end{equation}
with $\varepsilon>0$, the real-valued function $\mathcal{L}_R(\tilde{R})=\tr((I-\tilde{R})M) $ defined in (\ref{eqn:L_R}), and the real-valued function $\bar{\mathcal{L}}_p:\mathbb{R}^6 \to \mathbb{R}^+ $ defined as
\begin{equation}
\bar{\mathcal{L}}_p(\mathsf{x})  = \mathsf{x}\T P^{-1}  \mathsf{x}  \label{eqn:L_p2}.
\end{equation}
It is easy to verify that $\frac{1}{p_M} \|\mathsf{x}\|^2 \leq \bar{\mathcal{L}}_p \leq \frac{1}{p_m} \|\mathsf{x}\|^2$. Let $|x^c_2|_{\mathcal{A}_2}\geq 0$ denote the distance to the set $\mathcal{A}_2$ such that
$
|x^c_2|_{\mathcal{A}_2}^2  : = \textstyle \inf_{y=(\bar{X},I_3,0,\bar{t})\in \mathcal{A}_2}  (\|\bar{X}-\hat{X}\|_F^2 + |\tilde{R}|_I^2
+ \|\mathsf{x}\|^2 +\|\bar{t}-t\|^2 )
= \textstyle |\tilde{R}|_I^2 + \|\mathsf{x}\|^2 .
$ Recall (\ref{eqn:L_R_bound}), one has
\begin{align}
&\underline{\alpha} |x^c_2|_{\mathcal{A}_2}^2 \leq \mathcal{L}(x^c_2) \leq \bar{\alpha} |x^c_2|_{\mathcal{A}_2}^2    \label{eqn:L_x2}
\end{align}
where $\underline{\alpha}: = \min\{4\lambda_{m}^{\bar{M}},\frac{\varepsilon}{p_M}\}, \bar{\alpha}: =\max\{4\lambda_{M}^{\bar{M}},\frac{\varepsilon}{p_m}\}$.
Using the fact that  $\dot{P}^{-1} = -P^{-1} \dot{P}P^{-1}$, the time-derivative of $\bar{\mathcal{L}}_p$ in the flows is given by
\begin{align}
\dot{\bar{\mathcal{L}}}_p &=  \mathsf{x}\T  (P^{-1} A(t) + A(t)\T P^{-1} - 2 C\T Q(t)  C + \dot{P}^{-1}  )  \mathsf{x}  \nonumber \\
&~~~+ 2  \mathsf{x}\T P^{-1}   \nu  \nonumber\\
&\leq  - \mathsf{x}\T  P^{-1} V(t) P^{-1}   \mathsf{x} -\mathsf{x}\T  C\T Q(t)  C \mathsf{x} + 2 \mathsf{x}\T P^{-1}   \nu  \nonumber\\ 
&\leq - \frac{v_m}{p_M^2}  \mathsf{x} \T \mathsf{x}     +  \frac{4\sqrt{2} \|\g\| }{p_m}   \|\mathsf{x}\| |\tilde{R}|_I
\label{eqn:dot_L_p2}
\end{align}
where we made use of the facts $- \mathsf{x}\T   C\T Q(t)  C \mathsf{x} \leq 0$, $p_m I_6 \leq P \leq p_M I_6 $ and
$\|\nu\| \leq  \|I-\tilde{R}\|_F \|\mathsf{g}\|
 = 2\sqrt{2} \|\g\| |\tilde{R}|_I
$. In view of (\ref{eqn:dot_L_R}) and (\ref{eqn:dot_L_p2}), one has
\begin{align}
\dot{\mathcal{L}}(x^c_2) & \leq -\lambda_R |\tilde{R}|_I^2- \frac{\varepsilon v_m}{p_M^2}  \mathsf{x} \T \mathsf{x}     +  \frac{4\sqrt{2}\varepsilon \g }{p_m}   \|\mathsf{x}\| |\tilde{R}|_I   \nonumber \\
& = - \begin{bmatrix}
|\tilde{R}|_I & \|\mathsf{x}\|
\end{bmatrix} \underbrace{\begin{bmatrix}
	\lambda_R & -\frac{2\sqrt{2}\varepsilon \|\g\| }{p_m}\\
	-\frac{2\sqrt{2}\varepsilon \|\g\| }{p_m} & \frac{\varepsilon v_m}{p_M^2}
	\end{bmatrix}}_{P_4}  \begin{bmatrix}
|\tilde{R}|_I \\ \|\mathsf{x}\|
\end{bmatrix} \nonumber \\
&\leq -\lambda_F  \mathcal{L}(x^c_2),  \quad \forall x^c_2 \in \mathcal{F}^c_2  \label{eqn:dL_x2}
\end{align}
where $P_4$ is positive definite by choosing $\varepsilon<{\lambda_R v_m P_m^2}/{(8\|\g\|^2 p_M^2)}$, and $\lambda_F:= \lambda_{m}^{P_4}/\bar{\alpha}$ with $\bar{\alpha}$ given in (\ref{eqn:L_x2}).
Using the facts $\mathsf{x}^+ = \mathsf{x}$ and $\bar{\mathcal{L}}_p(\mathsf{x}^+)  =\bar{\mathcal{L}}_p (\mathsf{x}) $,  one can also show that
\begin{align}
\mathcal{L}({x^c_2}^+) - \mathcal{L}(x^c_2)
&=\mathcal{L}_R(\tilde{R}^+)  - \mathcal{L}_R(\tilde{R}) + \varepsilon \bar{\mathcal{L}}_p(\mathsf{x}^+)  - \varepsilon\bar{\mathcal{L}}_p (\mathsf{x}) \nonumber \\
&\leq -\delta, \quad \forall x^c_2 \in \mathcal{J}^c_2. \label{eqn:L_x+2}
\end{align}
Then, in view of (\ref{eqn:L_x2}), (\ref{eqn:dL_x2}) and (\ref{eqn:L_x+2}), the rest of the proof is completed using similar steps as in the proof of Theorem \ref{theo:theo_1}.

\subsection{Proof of Theorem \ref{theo:theo_3}} \label{sec:theo_3}
In view of (\ref{eqn:dX}), (\ref{eqn:Delta}),  (\ref{eqn:gamma}), (\ref{eqn:set Q}), and (\ref{eqn:observer_X3})-(\ref{eqn:innovation_term3}), one has the following hybrid closed-loop system:
\begin{equation}
\mathcal{H}^c_3:  \begin{cases}
\dot{x}^c_3 ~~= F_3(x^c_3),   & x^c_3\in \mathcal{F}^c_3 \\
{x^c_3}^+  = G_3(x^c_3),    &   x^c_3\in \mathcal{J}^c_3
\end{cases} \label{eqn:closed-loop3}
\end{equation}
where the flow and jump sets are defined as $\mathcal{F}^c_3:=\{ x^c_3=(x^c_1,\hat{b}_\omega, \tilde{b}_\omega ) \in \mathcal{S}^c_3: x^c_1 \in \mathcal{F}^c_1  \}$ and  $\mathcal{J}^c_3 :=\{ x^c_3=(x^c_1,\hat{b}_\omega, \tilde{b}_\omega) \in \mathcal{S}^c_3: x^c_1 \in \mathcal{J}^c_1 \}$, and the flow and jump maps are given by
\begin{align*}
&F_3(x^c_3)= \begin{pmatrix}
f(\hat{X},\omega_y-\hat{b}_\omega,a) - \Delta\hat{X}\\
 \tilde{R} ((\hat{R}\tilde{b}_\omega)^\times-k_R\mathbb{P}_a(M\tilde{R}))\\
 -k_pk_c\tilde{p}_e + \tilde{v}  -  (\tilde{R}\hat{R}\tilde{b}_\omega)^\times (p-p_c-\tilde{p}_e) \\
-  k_v k_c  \tilde{p}_e     - (\tilde{R}\hat{R}\tilde{b}_\omega)^\times (v-\tilde{v}) + (I_3-\tilde{R})\mathsf{g}\\
 1\\
-k_\omega   \hat{R}\T \psi(M\tilde{R}) \\
 -k_\omega   \hat{R}\T \psi(M\tilde{R})
\end{pmatrix}  \\
&G_3(x^c_3) = \begin{pmatrix}
(X_q^{-1} \hat{X})\T, ~
(\tilde{R}R_q)\T,~
\tilde{v}\T, ~
\tilde{p}_e\T, ~
t,~
\hat{b}_\omega\T,~
\tilde{b}_\omega\T
\end{pmatrix}\T
\end{align*}
where the following facts have been used: $\tilde{p}_e = (p-p_c)-\tilde{R}(\hat{p}-p_c)$, $\tilde{R}(\hat{R}\tilde{b}_\omega)^\times(\hat{p}-p_c)= (R\tilde{b}_\omega)^\times \tilde{R}(\hat{p}-p_c)= (R\tilde{b}_\omega)^\times (p-p_c - \tilde{p}_e)$, $\tilde{R}(\hat{R}\tilde{b}_\omega)^\times \hat{v} =(R\tilde{b}_\omega)^\times \tilde{R}\hat{v} = (R\tilde{b}_\omega)^\times (v-\tilde{v}) $.
Note that the sets $\mathcal{F}^c_3$ and $\mathcal{J}^c_3$ are closed, and $\mathcal{F}^c_3 \cup \mathcal{J}^c_3 = \mathcal{S}^c_3$. Note also that the closed-loop system (\ref{eqn:closed-loop3}) satisfies the hybrid basic conditions of \cite{goebel2009hybrid} and is autonomous by taking $\omega_y$, $a$, $p$ and $v$ as functions of time. By following similar steps as in the proof of Proposition \ref{pro:pro_1} and \ref{pro:undesired_eq}, one can show that all the undesired equilibrium points of the flows of $\mathcal{H}^c_3$ lie in the jump set $\mathcal{J}^c_3$, which is omitted here.

Consider the real-valued function $\mathcal{V}_R = \tr((I_3-\tilde{R})M) + \frac{1}{k_\omega} \tilde{b}_\omega\T \tilde{b}_\omega $, whose time-derivative in the flows is given by
\begin{align}
\dot{\mathcal{V}}_R  &=  \tr( -\tilde{R} ((\hat{R}\tilde{b}_\omega)^\times-k_R\mathbb{P}_a(M\tilde{R})) M)  -  2\tilde{b}_\omega\T  \hat{R}\T \psi(M\tilde{R})  \nonumber \\
 &= - k_R \|\mathbb{P}_a(M\tilde{R} )\|_F^2  \leq 0 \label{eqn:dotV_R}
\end{align}
where we made use of the facts   $\tr(-M\tilde{R}(\hat{R}\tilde{b}_\omega)^\times)=\langle\langle (\hat{R}\tilde{b}_\omega)^\times, M\tilde{R} \rangle\rangle = 2\psi(M\tilde{R})\T \hat{R} \tilde{b}_\omega$ and $ \tr(M\tilde{R}   \mathbb{P}_a(M\tilde{R})) =-\langle\langle  \mathbb{P}_a(M\tilde{R}), M\tilde{R} \rangle\rangle=-\langle\langle  \mathbb{P}_a(M\tilde{R}), \mathbb{P}_a(M\tilde{R} )\rangle\rangle$. Then, one concludes that $\mathcal{V}_R$ is non-increasing in the flows. By virtue of Proposition \ref{pro:undesired_eq}, for each jump one has
\begin{align}
\mathcal{V}_R^+ - \mathcal{V}_R  \leq  -\delta \label{eqn:V_R^+}.
\end{align}
 Therefore, for any $x^c_3(0,0)\in \mathcal{S}^c_3$, there exists  $c_{b_\omega} >0$ such that $c_{b_\omega}:= \sup_{(t,j)\succeq (0,0)} \|\tilde{b}_\omega(t,j)\|$ for all $(t,j)\in \dom x^c_3$. Note that $\|\tilde{b}_\omega(t,j)\|^2 \leq \mathcal{V}_R(t,j) \leq  \mathcal{V}_R(0,0)$, which implies that $c_{b_\omega}$ is bounded by the initial conditions. Let us modify the real-valued function $\bar{\mathcal{L}}_R:SO(3)\times \mathbb{R}^3 \to \mathbb{R}^+$ as follows:
\begin{equation}
\bar{\mathcal{L}}_R(\tilde{R},\tilde{b}_\omega)  =  \mathcal{L}_R(\tilde{R}) + \frac{1}{k_\omega} \tilde{b}_\omega\T \tilde{b}_\omega - \bar{\mu} \psi(\tilde{R})\T \hat{R}\tilde{b}_\omega  \label{eqn:L_R2}
\end{equation}
where $\bar{\mu}>0$. Let $e_1=[|\tilde{R}|_I,\|\tilde{b}_\omega\|]\T$. Following similar steps as in the proof of \cite[Theorem 1]{berkane2017hybrid} and \cite[Theorem 2]{wang2019hybrid}, there exists a constant $\bar{\mu}^*$ such that for all $\bar{\mu} \leq \bar{\mu}^*$ one has
\begin{align}
&\underline{c}_R \|e_1\|^2 \leq \bar{\mathcal{L}}_R \leq \bar{c}_R \|e_1\|^2 \label{eqn:L_R_bound2}, \\
&\dot{\bar{\mathcal{L}}}_R \leq -  \bar{\lambda}_{R} \|e_1\|^2       \qquad  \forall x^c_3 \in \mathcal{F}^c_3 \label{eqn:dot_L_R2}
\end{align}
for some positive constants $\underline{c}_R,\bar{c}_R $ and $ \bar{\lambda}_{R}$. From \cite[Theorem 1]{berkane2017hybrid} and \cite[Theorem 2]{wang2019hybrid}, the constant $\bar{\lambda}_{R}$ depends on $c_{b_\omega}$, which is associated to the initial conditions.

On the other hand, we consider the real-valued functions $\mathcal{L}_p$ defined in (\ref{eqn:Lp}). Defining $e_2 := [\|\tilde{p}_e\|~ \|\tilde{v}\|]\T$, one verifies that $e_2\T P_1 e_2 \leq \mathcal{L}_p \leq e_2\T P_2 e_2$ as shown in (\ref{eqn:L_p_bound}). From Assumption \ref{assum:2}, there exist two constants $c_p, c_v$ such that $c_p := \sup_{t\geq 0} \|p-p_c\|, c_v := \sup_{t\geq 0} \|v\| $. Then, in the flows of (\ref{eqn:closed-loop3}) one has
\begin{align*}
\frac{d}{dt} \|\tilde{p}_e\|^2  & =  2\tilde{p}_e\T(-k_c k_p \tilde{p}_e + \tilde{v}-  (R\tilde{b}_\omega)^\times (p-p_c-\tilde{p}_e)) \\
& \leq  -2k_c k_p \|\tilde{p}_e\|^2 + 2c_p \|\tilde{p}_e\| \|\tilde{b}_\omega \| + 2\tilde{p}_e\T \tilde{v} \\
\frac{d}{dt}  \|\tilde{v}\|^2  & =    2 \tilde{v}\T (  - k_c k_v  \tilde{p}_e   -  (R\tilde{b}_\omega)^\times (v-\tilde{v}) + (I-\tilde{R})\mathsf{g})  \nonumber \\
& \leq 2(c_v \|\tilde{b}_\omega\| + 2\sqrt{2}\|\g\| |\tilde{R}|_I)\|\tilde{v}\| - 2k_c k_v \tilde{v}\T\tilde{p}_e \\ %
-\frac{d}{dt}  \tilde{p}_e\T \tilde{v} & =  (k_c k_p \tilde{p}_e - \tilde{v}+  (R\tilde{b}_\omega)^\times (p-p_c-\tilde{p}_e))\T \tilde{v} \\
& ~~ + \tilde{p}_e\T  (  k_c k_v \tilde{p}_e   +  (R\tilde{b}_\omega)^\times (v-\tilde{v})   - (I_3-\tilde{R})\mathsf{g}) \\
& \leq - \|\tilde{v}\|^2 +  k_c k_p \|\tilde{p}_e\|\tilde{v}\| + c_p \|\tilde{v}\|\|\tilde{b}_\omega\|   \\
& ~~ + k_c k_v \|\tilde{p}_e\|^2 + c_v \|\tilde{b}_\omega\|\|\tilde{p}_e\| +2\sqrt{2}\|\g\||\tilde{R}|_I \|\tilde{p}_e\|
\end{align*}
where  we made use of the facts: $\|I_3-\tilde{R}\|_F = 2\sqrt{2}|\tilde{R}|_I$, $((R\tilde{b}_\omega)^\times)\T = - (R\tilde{b}_\omega)^\times$, $x\T (R\tilde{b}_\omega)^\times x =0,  \forall x\in \mathbb{R}^3$. Let $c_2: =  \max\{c_p+\mu c_v,\frac{c_v}{k_ck_v} + \mu c_p\}, c_3:= 2\sqrt{2}\|\g\|\max\{\frac{1}{k_ck_v},\mu  \}$ and $c_4 := \max\{c_2,c_3\}$. Then,  the time-derivative of $\mathcal{L}_p$ in the flows of (\ref{eqn:closed-loop3}) satisfies
\begin{align}
\dot{\mathcal{L}}_p
& \leq    -(k_p-\mu k_v ) k_c \|\tilde{p}_e\|^2  - \mu \|\tilde{v}\|^2   + \mu k_c k_p \|\tilde{p}_e\| \|\tilde{v}\| \nonumber \\
& \quad +  c_2 ( \|\tilde{p}_e\|   +   \|\tilde{v}\| )\|\tilde{b}_\omega\| + c_3  ( \|\tilde{v}\|+   \|\tilde{p}_e\|) |\tilde{R}|_I \nonumber \\
& \leq  - e_2\T P_3 e_2 + 2c_4 \|e_1\| \|e_2\|  \label{eqn:dot_L_p3}
\end{align}
where $P_3$ is given in (\ref{eqn:dot_L_p}), and the following facts have been used: $|\tilde{R}|_I + \|\tilde{b}_\omega\| \leq   \sqrt{2}\|e_1\|$ and $\|\tilde{v}\| + \|\tilde{p}_e\| \leq  \sqrt{2}\|e_2\|$. Pick
$$0< \mu < \min\left\{\frac{1}{\sqrt{k_c k_v}}, \frac{4k_p}{4k_v + k_ck_p^2} \right\}$$
such that the matrices $P_1,P_2$ and $P_3$ are positive definite.

Consider the following Lyapunov function candidate:
\begin{equation}
\mathcal{L}(x^c_3): = \bar{\mathcal{L}}_R(\tilde{R},\tilde{b}_\omega) + \varepsilon\mathcal{L}_p(\tilde{p}_e,\tilde{v})  \label{eqn:Lyp3}
\end{equation}
where $\varepsilon>0$.
 From (\ref{eqn:L_R_bound2}) and (\ref{eqn:Lp}), one has
\begin{align}
&\underline{\alpha} \|x^c_3\|_{\mathcal{A}_3}^2 \leq \mathcal{L} (x^c_3) \leq \bar{\alpha}  \|x^c_3\|_{\mathcal{A}_3}^2   \label{eqn:L_x3}
\end{align}
where  $\underline{\alpha}: = \min\{\underline{c}_R,\varepsilon \lambda_{m}^{P_1}\}, \bar{\alpha}: =\max\{\bar{c}_R,\varepsilon\lambda_{M}^{P_2}\}$. From (\ref{eqn:dot_L_R2}) and (\ref{eqn:dot_L_p3}), for any $x^c_3\in \mathcal{F}^c_3$ one has
\begin{align}
\dot{\mathcal{L}}(x^c_3)  & \leq -\bar{\lambda}_{R} \|e_1\|^2 - \varepsilon\lambda_{m}^{P_3}\|e_2\|^2 + 2\varepsilon c_4 \|e_1\|\|e_2\|  \nonumber \\
& = - \begin{bmatrix}
\|e_1\| & \|e_2\|
\end{bmatrix} \underbrace{\begin{bmatrix}
	\bar{\lambda}_{R} & -\varepsilon c_4\\
	-\varepsilon c_4 & \varepsilon\lambda_{m}^{P_3}
	\end{bmatrix}}_{P_4}  \begin{bmatrix}
\|e_1\| \\ \|e_2\|
\end{bmatrix} \nonumber \\
&\leq -\lambda_F  \mathcal{L}(x^c_3)  \label{eqn:dL_x3}
\end{align}
where $P_4$ is positive definite by choosing $\varepsilon<{\bar{\lambda}_{R}\lambda_{m}^{P_3}}/{c_4^2}$, and $\lambda_F:= \lambda_{m}^{P_4}/\bar{\alpha}$.  In view of  (\ref{eqn: definition_mu_R})-(\ref{eqn:J_map}), (\ref{eqn:Lp}), (\ref{eqn:closed-loop3}) and (\ref{eqn:L_R2}) , for any $x^c_3\in \mathcal{J}^c_3$ one has	
 \begin{align*}
& \mathcal{L}({x^c_3}^+) - \mathcal{L}(x^c_3) \nonumber \\
& = \bar{\mathcal{L}}_R(\tilde{R}^+,\tilde{b}_\omega^+)  - \bar{\mathcal{L}}_R(\tilde{R},\tilde{b}_\omega) + \varepsilon \mathcal{L}_p(\tilde{p}_e^+,\tilde{v}^+)  - \varepsilon\mathcal{L}_p(\tilde{p}_e,\tilde{v}) \\
&=-\delta - \bar{\mu} \psi(\tilde{R})\T\hat{R} \tilde{b}_\omega + \bar{\mu} \psi(\tilde{R}R_q)\T R_q\T\hat{R} \tilde{b}_\omega   \\
&\leq  -\delta + 4 \bar{\mu} c_{b_\omega}
 \end{align*}
where we made use of the results from (\ref{eqn:L_x+}), and the fact $\|\psi(\tilde{R})\| \leq 1$ from $\|\psi(\tilde{R})\|^2 = 4|\tilde{R}|_I^2|(1-|\tilde{R}|_I^2)\leq 1, \forall |\tilde{R}|_I^2\in[0,1]$. Choosing $\bar{\mu}< \min\{\bar{u}^*, \delta/2c_{b_\omega}\}$, one has
 \begin{align}
 \mathcal{L}({x^c_3}^+) - \mathcal{L}(x^c_3)  \leq - \delta^*, \quad  \forall x^c_3\in \mathcal{J}^c_3 \label{eqn:L_x+3}
 \end{align}
where $\delta^*:=-\delta + 2 \bar{\mu} c_{b_\omega} >0$. In view of (\ref{eqn:dL_x3}) and (\ref{eqn:L_x+3}), one has $0\leq \mathcal{L}(x^c_3(t,j))\leq   \mathcal{L}(x^c_3(0,0))-j\delta^* $, which leads to $j\leq J := \left\lceil  {\mathcal{L}(x^c_3(0,0))}/{\delta^*} \right\rceil$. This implies that the number of jumps is finite. Moreover, one has $\mathcal{L}(x^c_3(t,j))\leq  \exp(-\lambda_F t)  \mathcal{L}(x^c_3(0,0)) \leq \exp (\lambda_F J) \exp(-\lambda_F (t+j))  \mathcal{L}(x^c_3(0,0)) $. Substituting (\ref{eqn:L_x3}), one concludes that for each $(t,j)\in \dom x^c_3$,
\begin{equation*}
|x^c_3(t,j)|^2_{\mathcal{A}_3} \leq \kappa\exp\left(- \lambda_F(t+j)\right) |x^c_3(0,0)|^2_{\mathcal{ A}_3}
\end{equation*}
where $ \kappa:=\exp( \lambda_F J) \bar{\alpha}/ \underline{\alpha} $. This completes the proof.

\subsection{Proof of Theorem \ref{theo:theo_4}} \label{sec:theo4}
The proof of Theorem \ref{theo:theo_4} is similar to the proof of Theorem \ref{theo:theo_2} and Theorem \ref{theo:theo_3}.
In view of (\ref{eqn:CRE}),  (\ref{eqn:gamma}), (\ref{eqn:set Q}), and (\ref{eqn:observer_X4})-(\ref{eqn:dx2}), one obtains the following hybrid closed-loop system:
\begin{equation}
\mathcal{H}^c_4:  \begin{cases}
\dot{x}^c_4 ~~= F_4(x^c_4),   & x^c_4\in \mathcal{F}^c_4 \\
{x^c_4}^+  = G_4(x^c_4),    &   x^c_4\in \mathcal{J}^c_4
\end{cases} \label{eqn:closed-loop4}
\end{equation}
where the flow and jump sets are defined as $\mathcal{F}^c_4:=\{ x^c_4=(x^c_2,\hat{b}_\omega,\tilde{b}_\omega) \in \mathcal{S}^c_4: x^c_2 \in \mathcal{F}^c_2 \}$ and  $\mathcal{J}^c_4 :=\{x^c_4=(x^c_2,\hat{b}_\omega,\tilde{b}_\omega) \in \mathcal{S}^c_4: x^c_2 \in \mathcal{J}^c_2\}$ with $\mathcal{F}^c_2,\mathcal{J}^c_2$ given in (\ref{eqn:closed-loop2}), and the flow and jump maps are given by
\begin{align*}\small 
F_4(x^c_4)= \begin{pmatrix}
f(\hat{X},\omega_y-\hat{b}_\omega,a) - \Delta\hat{X}\\
\tilde{R} (-k_R\mathbb{P}_a(M\tilde{R}))\\
A(t) \mathsf{x} - L C \mathsf{x} + \nu \\
1\\
-k_\omega   \hat{R}\T \psi(M\tilde{R}) \\
-k_\omega   \hat{R}\T \psi(M\tilde{R})
\end{pmatrix},
G_4(x^c_4) = \begin{pmatrix}
X_q^{-1} \hat{X}\\
\tilde{R}R_q\\
\mathsf{x}\\
t \\
\hat{b}_\omega\\
\tilde{b}_\omega
\end{pmatrix}.
\end{align*}
Note that the sets $\mathcal{F}^c_4$ and $\mathcal{J}^c_4$ are closed, and $\mathcal{F}^c_4 \cup \mathcal{J}^c_4 = \mathcal{S}^c_4$. Note also that the closed-loop system (\ref{eqn:closed-loop2}) satisfies the hybrid basic conditions of \cite{goebel2009hybrid} and is autonomous by taking $\omega_y$, $a$, $A(t)$ and $L$ as functions of $t$. By following similar steps as in the proof of Proposition \ref{pro:pro_1} and \ref{pro:undesired_eq}, one can show that all the undesired equilibrium points of the flows of $\mathcal{H}^c_4$ lie in the jump set $\mathcal{J}^c_4$, which is omitted here.

Consider the following Lyapunov function candidate:
\begin{equation}
\mathcal{L}(x^c_4): = \bar{\mathcal{L}}_R(\tilde{R},\tilde{b}_\omega) + \varepsilon\bar{\mathcal{L}}_p(\mathsf{x})  \label{eqn:Lyp4}
\end{equation}
where $\varepsilon>0$, the real-valued function $\bar{\mathcal{L}}_R$ is defined in (\ref{eqn:L_R2}) and the real-valued function $\bar{\mathcal{L}}_p$ is defined in (\ref{eqn:L_p2}).
It is easy to verify that $\frac{1}{p_M} \|\mathsf{x}\|^2 \leq \bar{\mathcal{L}}_p \leq \frac{1}{p_m} \|\mathsf{x}\|^2$. Using the fact $\frac{1}{p_M} \|\mathsf{x}\|^2\leq \bar{\mathcal{L}}_p \leq \frac{1}{p_m} \|\mathsf{x}\|^2$ and property (\ref{eqn:L_R_bound2}), one has
\begin{align}
&\underline{\alpha} |x^c_4|_{\mathcal{A}_4}^2 \leq \mathcal{L}(x^c_4) \leq \bar{\alpha} |x^c_4|_{\mathcal{A}_4}^2  \label{eqn:L_x4}
\end{align}
where $\underline{\alpha}: = \min\{\underline{c}_R,\frac{\varepsilon}{p_M}\}, \bar{\alpha}: =\max\{\bar{c}_R,\frac{\varepsilon}{p_m}\}$. Using the facts $c_p := \sup_{t\geq 0} \|p-p_c\|$, $ c_v := \sup_{t\geq 0} \|v\| $ and $\|e_1\|^2=|\tilde{R}|_I^2 + \|\tilde{b}_\omega\|^2 $, one can show that $\|\nu\|^2 \leq (c_p + c_v) \|\tilde{b}_\omega\|^2 + 8\|\mathsf{g}\|^2 |\tilde{R}|_I^2 \leq c_5^2 \|e_1\|^2$ with $c_{5}:=\max\{\sqrt{c_p + c_v},2\sqrt{2}\|\mathsf{g}\|\}$.
Then, the time-derivative of $\bar{\mathcal{L}}_p$ in the flows is given by
\begin{align}
\dot{\bar{\mathcal{L}}}_p &=  \mathsf{x}\T  (P^{-1} A(t) + A(t)\T P^{-1} - 2 C\T Q(t)  C + \dot{P}^{-1} )  \mathsf{x}  \nonumber \\
& ~~~+ 2  \mathsf{x}\T P^{-1}   \nu \nonumber \\
&\leq  - \mathsf{x}\T  P^{-1} V(t) P^{-1}   \mathsf{x}  -\mathsf{x}\T  C\T Q(t)  C \mathsf{x} + 2 \mathsf{x}\T P^{-1}   \nu  \nonumber \\
&\leq -   \frac{v_m }{p_M^2}    \mathsf{x} \T \mathsf{x}     +  \frac{2c_5 }{p_m}   \|\mathsf{x}\| \|e_1\|
\label{eqn:dot_L_p4}
\end{align}
where we made use of the facts   $- \mathsf{x}\T   C\T Q(t)  C \mathsf{x} \leq 0$ and $\|\nu\| \leq c_{5} \|e_1\|$. From (\ref{eqn:dot_L_R2}) and (\ref{eqn:dot_L_p4}), one obtains
\begin{align}
\dot{\mathcal{L}}(x^c_4) & \leq -\bar{\lambda}_{R} \|e_1\|^2- \frac{\varepsilon v_m}{p_M^2}  \mathsf{x} \T \mathsf{x}     +  \frac{2 \varepsilon c_5 }{p_m} \|e_1\|  \|\mathsf{x}\|     \nonumber \\
& = - \begin{bmatrix}
\|e_1\| & \|\mathsf{x}\|
\end{bmatrix} \underbrace{\begin{bmatrix}
	\bar{\lambda}_{R} & -\frac{ \varepsilon c_5 }{p_m}\\
	-\frac{  \varepsilon c_5 }{p_m} & \frac{\varepsilon v_m}{p_M^2}
	\end{bmatrix}}_{P_4}  \begin{bmatrix}
\|e_1\| \\ \|\mathsf{x}\|
\end{bmatrix} \nonumber \\
&\leq -\lambda_F  \mathcal{L}(x^c_4),  \quad  \forall x^c_4 \in \mathcal{F}^c_4  \label{eqn:dL_x4}
\end{align}
where $P_4$ is positive definite by choosing $\varepsilon<{\bar{\lambda}_{R} v_m p_m^2}/{( c_5^2 p_M^2)}$, and $\lambda_F:= \lambda_{m}^{P_4}/\bar{\alpha}$ with $\bar{\alpha}$ given in (\ref{eqn:L_x4}). In view of  (\ref{eqn: definition_mu_R})-(\ref{eqn:J_map}), (\ref{eqn:L_p2}), (\ref{eqn:L_R2}) and (\ref{eqn:closed-loop4}), for any $x^c_4\in \mathcal{J}^c_4$ one has	
\begin{align*}
&\mathcal{L}({x^c_4}^+) - \mathcal{L}(x^c_4) \\
& = \bar{\mathcal{L}}_R(\tilde{R}^+,\tilde{b}_\omega^+)  - \bar{\mathcal{L}}_R(\tilde{R},\tilde{b}_\omega) + \varepsilon \bar{\mathcal{L}}_p(\mathsf{x}^+)  - \varepsilon\bar{\mathcal{L}}_p(\mathsf{x}) \\
&=-\delta - \bar{\mu} \psi(\tilde{R})\T\hat{R} \tilde{b}_\omega + \bar{\mu} \psi(\tilde{R}R_q)\T R_q\T\hat{R} \tilde{b}_\omega   \\
&\leq  -\delta + 2 \bar{\mu} c_{b_\omega}
\end{align*}
where we made use of the results from (\ref{eqn:L_x+}),  and the fact $\|\psi(\tilde{R})\| \leq 1, \forall \tilde{R}\in SO(3)$. Choosing $\bar{\mu}< \min\{\bar{u}^*, \delta/2c_{b_\omega}\}$, one has
\begin{align}
\mathcal{L}({x^c_4}^+) - \mathcal{L}(x^c_4)  \leq - \delta^*, \quad  \forall x^c_4\in \mathcal{J}^c_4 \label{eqn:L_x+4}
\end{align}
where $\delta^*:=-\delta + 2 \bar{\mu} c_{b_\omega} >0$. In view of (\ref{eqn:L_x4}), (\ref{eqn:dL_x4}) and (\ref{eqn:L_x+4}), the rest of the proof can be completed by using similar steps as in the proof of Theorem \ref{theo:theo_3}.

\subsection{Proof of Lemma \ref{lem:existenceP3}}\label{sec:existenceP3}
For the sake of simplicity, let $\varpi: = \omega_y - \hat{b}_\omega = \omega - \tilde{b}_\omega$.
From (\ref{eqn:closed-loop_F5}) and (\ref{eqn:dotV_R})-(\ref{eqn:V_R^+}), one shows that $\tilde{b}_\omega$ and  $\dot{\tilde{b}}_\omega$ are bounded. Moreover, from Assumption \ref{assum:2} and   \ref{assum:3}, one obtains that $\omega$ and $\dot{\omega}$ are uniformly bounded.
Hence, one can show that the derivations of $\varpi$ and $A(t)$ are well-defined and bounded for all $t\geq 0$. Define $N_0 = C$, $N_1 = N_0A(t)$ and $N_2 = N_1 A(t) + \dot{N}_1$.
From \cite[Theorem 4]{bristeau2010design}, the pair $(A(t),C)$ is uniformly observable   if there exists a  positive constant $\bar{\mu}$ such that
\begin{equation}
  \mathcal{O}(\tau)\T \mathcal{O}(\tau) d\tau \geq \bar{\mu} I_9, \quad \forall t\geq 0 \label{eqn:conditonO}
\end{equation}
where
\begin{align*}
\mathcal{O}(t) =  \begin{bmatrix}
N_0 \\
N_1 \\
N_2
\end{bmatrix} = \begin{bmatrix}
I_3 & 0_{3\times 3} & 0_{3\times 3}  \\
- \varpi^\times & I_3 & 0_{3\times 3}  \\
(\varpi^\times)^2 - \dot{\varpi}^\times &  - 2\varpi^\times  & I_3
\end{bmatrix}.
\end{align*}
It is easy to verify that matrix $\mathcal{O}$ is well-defined and has full rank since $\det(\mathcal{O}) =1 $ for all $t\geq 0$.  This implies that there exists a constant $\bar{\mu}>0$ such that $\mathcal{O}(t)\T \mathcal{O}(t) \geq \bar{\mu} I_9, \forall t\geq 0$. Therefore,  one can conclude that the pair $(A(t),C)$ is uniformly observable. This completes the proof.

\bibliographystyle{ieeetran}
\bibliography{mybib}

\begin{IEEEbiography}[{\includegraphics[width=1in,height=1.25in,clip,keepaspectratio]{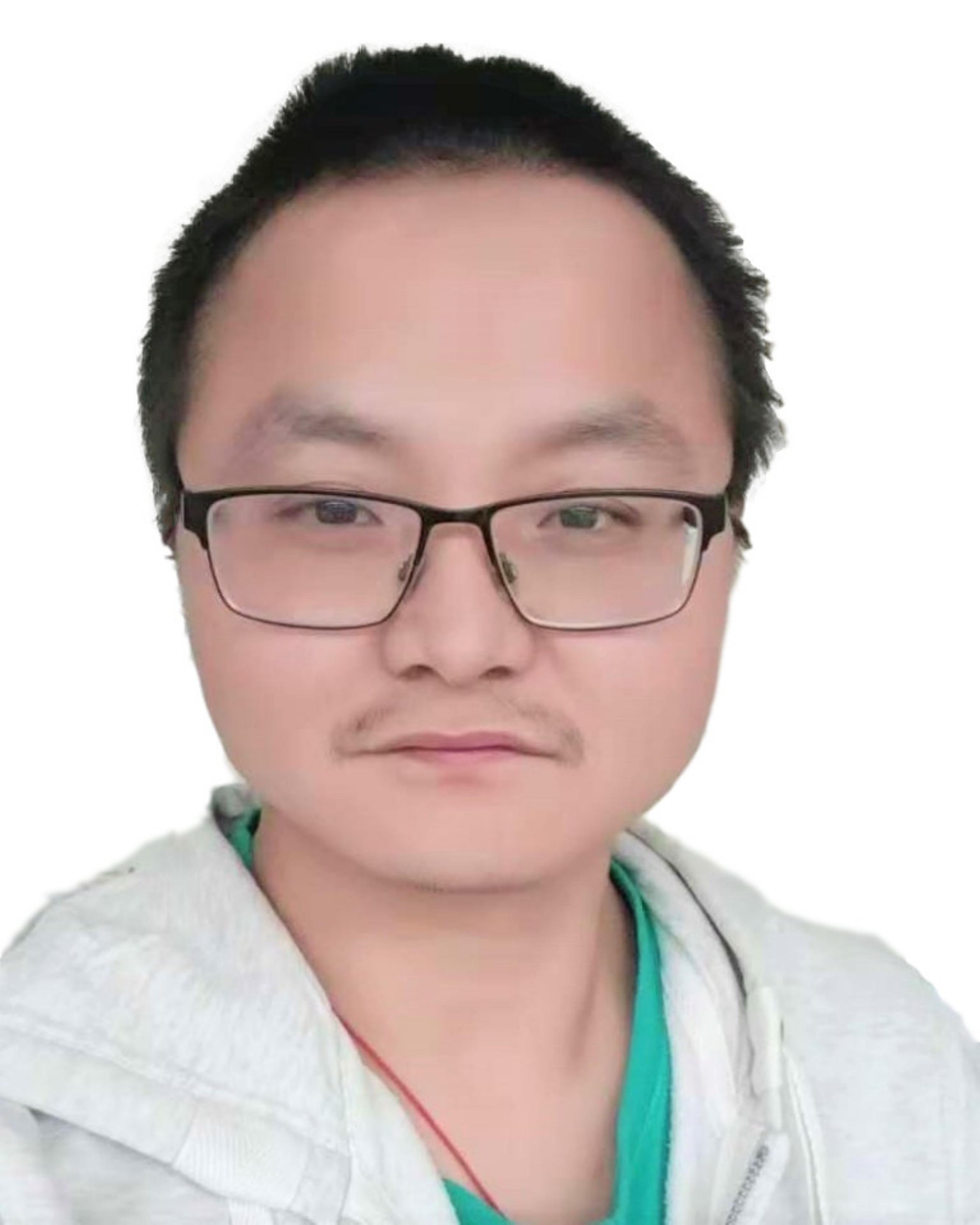}}]{Miaomiao Wang} received his B.Sc.   in control science and engineering from Huazhong University of Science and Technology, China, in 2013, and his M.Sc.   in Control Engineering from Lakehead University, Canada, in 2015.
He is currently a Ph.D student and a research assistant in the department of Electrical and Computer Engineering at Western University, Canada. He has received the prestigious Ontario Graduate Scholarship (OGS) at Western University in 2018. His current  research interests are in the areas of geometric control and nonlinear estimation with applications to autonomous robotic systems.
\end{IEEEbiography}

\begin{IEEEbiography}[{\includegraphics[width=1in,height=1.25in,clip,keepaspectratio]{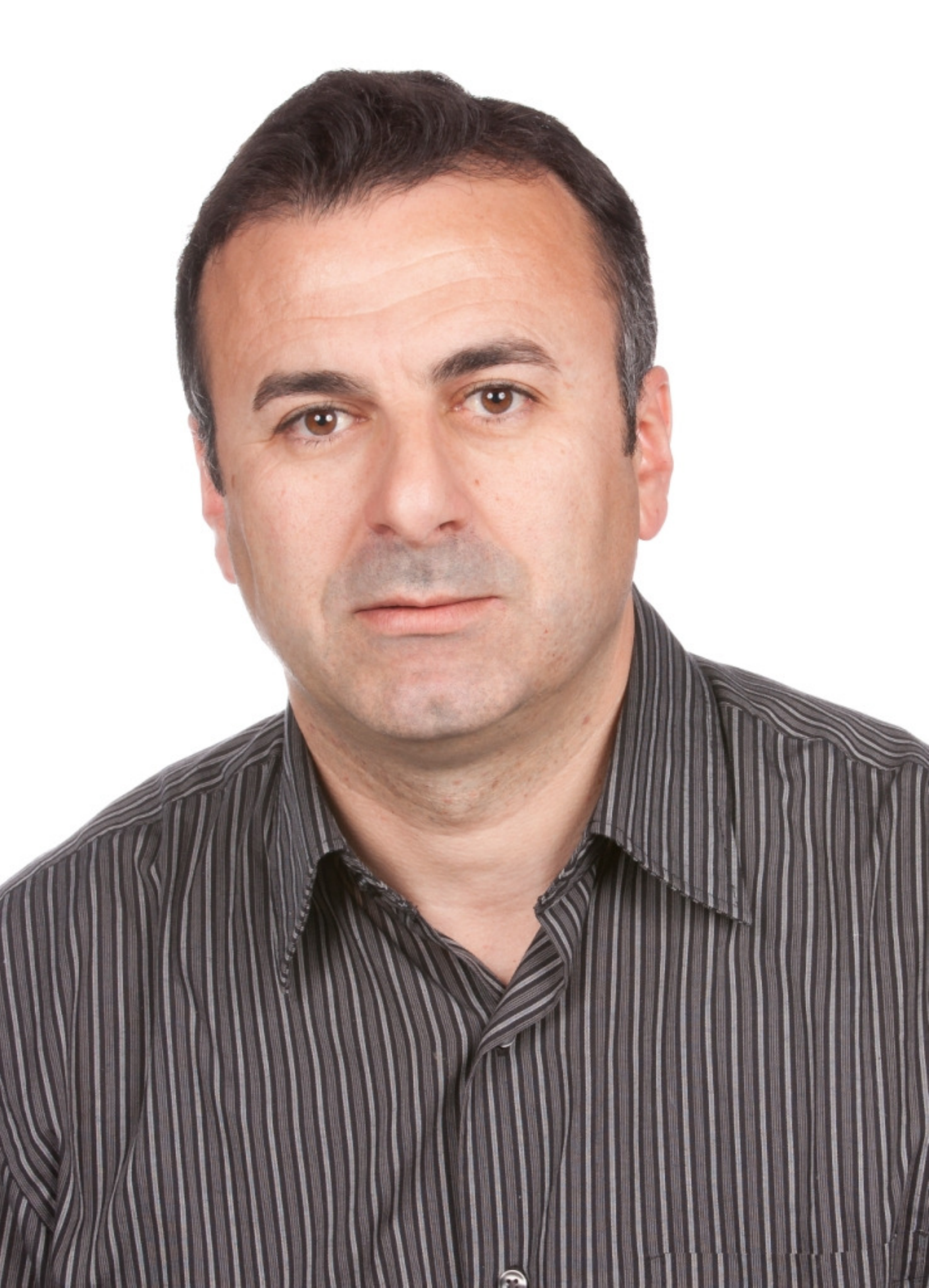}}]{Abdelhamid Tayebi}
received his B. Sc. in Electrical Engineering from Ecole Nationale Polytechnique, Algiers, in 1992, his M. Sc. (DEA) in robotics from Universit\'e Pierre \& Marie Curie, Paris, France in 1993, and his Ph. D. in Robotics and Automatic Control from Universit\'e de Picardie Jules Verne, France in 1997. He joined the department of Electrical Engineering at Lakehead University in 1999 where he is presently a Professor.
He is a Senior Member of IEEE and serves as an Associate Editor for Automatica and IEEE Transactions on Control Systems Technology. He is the founder and Director of the Automatic Control Laboratory at Lakehead University. His current research interests are in the broad area of Control Systems, Cooperative Control, Iterative Learning Control and Aerial Robotics.
\end{IEEEbiography}


\end{document}